\documentclass[aop,MSNbibl,seceqn,dvips]{arximspdf}

\doi{10.1214/12-AOP782} %kopijuoti is PTS
\volume{41}
\issue{4}
\pubyear{2013}
\firstpage{2755}
\lastpage{2790}

\makeatletter
\newcommand{\rrvert}{\vert}
\newcommand{\llvert}{\vert}

\newtheorem{lemma}{Lemma}[section]
\newtheorem{theorem}[lemma]{Theorem}
\newtheorem{proposition}[lemma]{Proposition}

\newproclaim{remark}[lemma]{Remark}

\def\S{{\mathcal{S}}}

\renewcommand{\a}{\alpha}
\renewcommand{\b}{\beta}
\renewcommand{\d}{\delta}
\newcommand{\g}{\gamma}
\newcommand{\eps}{\varepsilon}
\newcommand{\eb}{{\mathbb E}}
\newcommand{\pb}{{\mathbb P}}

\newcommand{\one}{{\mathbf1}}

\newcommand{\wt}{\widetilde}

\newcommand{\var}{\operatorname{var}}
\newcommand{\stas}{\stackrel{\mathrm{a.s.}}{\rightarrow}}
\newcommand{\eqd}{\stackrel{d}{=}}

\newcommand{\nto}{n\to\infty}
\newcommand{\xto}{x\to\infty}
\newcommand{\uto}{u\to\infty}
\newcommand{\vep}{\varepsilon}

\newcommand{\bbr}{{\mathbb R}}
\newcommand{\bbz}{{\mathbb Z}}

\renewcommand{\P}{\pb}
\newcommand{\E}{\eb}

\makeatother

\begin{document}
\begin{frontmatter}

\title{Large deviations for solutions to stochastic
recurrence equations under Kesten's condition\thanksref{T1}}
\runtitle{Large deviations}

\thankstext{T1}{Supported by NCN Grant UMO-2011/01/M/ST1/04604.}

\begin{aug}
\author[A]{\fnms{D.} \snm{Buraczewski}\corref{}\ead[label=e1]{dbura@math.uni.wroc.pl}},
\author[A]{\fnms{E.} \snm{Damek}\ead[label=e2]{edamek@math.uni.wroc.pl}},
\author[B]{\fnms{T.} \snm{Mikosch}\thanksref{t2}\ead[label=e3]{mikosch@math.ku.dk}}
\and
\author[A]{\fnms{J.} \snm{Zienkiewicz}\ead[label=e4]{zenek@math.uni.wroc.pl}}
\runauthor{Buraczewski, Damek, Mikosch and Zienkiewicz}
\affiliation{Uniwersytet Wroclawski, Uniwersytet Wroclawski,
University~of~Copenhagen  and Uniwersytet Wroclawski}
\address[A]{D. Buraczewski\\
E. Damek\\
J. Zienkiewicz\\
Instytut Matematyczny\\
Uniwersytet Wroclawski\\
50-384 Wroclaw\\
pl. Grunwaldzki 2/4\\
Poland\\
\printead{e1}\\
\hphantom{E-mail: }\printead*{e2}\\
\hphantom{E-mail: }\printead*{e4}} %adresu isvedimo komanda gale!
\address[B]{T. Mikosch\\
University of Copenhagen\\
Universitetsparken 5\\
DK-2100 Copenhagen\\
Denmark\\
\printead{e3}}
\end{aug}

\thankstext{t2}{Supported in part by the Danish
Natural Science Research Council (FNU) Grants
09-072331 ``Point process modeling and statistical inference''
and 10-084172 ``Heavy tail phenomena: Modeling and estimation.''}

% HISTORY:
\received{\smonth{5} \syear{2011}}
\revised{\smonth{1} \syear{2012}}

% ABSTRACT
%
\begin{abstract}
In this paper we prove large deviations results for partial sums
constructed from the solution to a stochastic recurrence equation. We
assume Kesten's condition [\textit{Acta Math.} \textbf{131} (1973)
207--248] under which the solution of the stochastic recurrence
equation has a marginal distribution with power law tails, while the
noise sequence of the equations can have light tails. The results of
the paper are analogs to those obtained by A. V. Nagaev [\textit{Theory
Probab. Appl.} \textbf{14} (1969) 51--64; 193--208] and S. V. Nagaev
[\textit{Ann. Probab.} \textbf{7} (1979) 745--789] in the case of
partial sums of i.i.d. random variables. In the latter case, the large
deviation probabilities of the partial sums are essentially determined
by the largest step size of the partial sum. For the solution to a
stochastic recurrence equation, the magnitude of the large deviation
probabilities is again given by the tail of the maximum summand, but
the exact asymptotic tail behavior is also influenced by clusters of
extreme values, due to dependencies in the sequence. We apply the large
deviation results to study the asymptotic behavior of the ruin
probabilities in the model.
\end{abstract}

% KEYWORDS

\end{frontmatter}

%s1 #&#
\section{Introduction}
Throughout the last 40 years, the stochastic recurrence equation
%
%e1.1 #&#
\begin{equation}
\label{eqkesten} Y_{n}= A_{n} Y_{n-1}+
B_{n},\qquad n \in\bbz,
\end{equation}
and its stationary solution have attracted much attention. Here
$(A_i,B_i)$, $i\in\bbz$, is an i.i.d. sequence, $A_i>0$ a.s., and
$B_i$ assumes real values.
[In what follows, we write $A,B,Y,\ldots\,$, for generic
elements of the strictly stationary sequences $(A_i)$, $(B_i)$,
$(Y_i),\ldots\,$, and we also write $c$ for any positive constant whose
value is not of interest.]

It is well known that if $\eb\log A <0$ and ${\eb\log^+}|B|<\infty$,
there exists a unique, strictly stationary
ergodic solution $(Y_i)$ to the stochastic recurrence equation (\ref
{eqkesten}) with representation
\[
Y_n=\sum_{i=-\infty}^n
A_{i+1}\cdots A_n B_i,\qquad n\in\bbz,
\]
where, as usual, we interpret the summand for $i=n$ as $B_n$.

One of the most interesting results for the stationary solution
$(Y_i)$ to the stochastic recurrence equation (\ref{eqkesten}) was
discovered by Kesten
\cite{kesten1973}. He proved under general conditions
that the marginal distributions of $(Y_i)$ have
power law tails. For later use, we formulate a version of this result
due to Goldie~\cite{goldie1991}.
%
%th1.1 #&#
\begin{theorem}[(Kesten~\cite{kesten1973},
Goldie~\cite{goldie1991})]\label{thmkesten}  Assume that the following
conditions hold:
\begin{itemize}
\item There exists $\a>0$ such that
%
%e1.2 #&#
\begin{equation}
\label{eqkestena} \eb A^\alpha=1.
\end{equation}
\item$\rho= \eb(A^\a\log A)$ and $\eb|B|^\a$ are both finite.
\item The law of $\log A$ is nonarithmetic.
\item For every $x$, $\pb\{Ax+B=x\}<1$.
\end{itemize}
Then $Y$ is regularly varying with index
$\alpha>0$. In particular, there exist constants
$c_\infty^+,c_\infty^-\ge0$ such that $c_\infty^++c_\infty^->0$ and
%
%e1.3 #&#
\begin{equation}
\label{eqpower} \pb\{Y>x\}\sim c_\infty^+ x^{-\alpha}\quad
\mbox{and}\quad
\pb\{Y\le-x\}\sim c_\infty^- x^{-\alpha} \qquad\mbox{as } \xto.
\end{equation}
Moreover, if $B\equiv1$ a.s., then the constant $c_\infty^+$ takes on
the form
\[
c_\infty:=\eb\bigl[(1+Y)^\alpha-Y^\alpha\bigr]/(\a
\rho).
\]
\end{theorem}

Goldie
\cite{goldie1991} also showed that similar results remain valid
for the stationary solution to stochastic recurrence equations of the
type $Y_{n}=
f(Y_{n-1},A_{n},B_{n}) $ for suitable functions $f$ satisfying some
contraction condition.

The power law tails (\ref{eqpower}) stimulated research on the
extremes of the sequence~$(Y_i)$. Indeed, if $(Y_i)$ were i.i.d. with tail
(\ref{eqpower}) and $c_\infty^+>0$, then the maximum sequence
$M_n=\max(Y_1,\ldots,Y_n)$ would
satisfy the limit relation
%
%e1.4 #&#
\begin{equation}
\label{eqextreme} \lim_{\nto}\pb\bigl\{ \bigl(c_\infty^+ n
\bigr)^{-1/\alpha} M_n\le x\bigr\}= e^{-x^{-\alpha}}=
\Phi_\alpha(x),\qquad x>0,
\end{equation}
where $\Phi_\alpha$ denotes the Fr\'echet distribution, that is, one of the
classical extreme value distributions; see Gnedenko~\cite{gnedenko1943};
cf. Embrechts et al.~\cite{embrechtskluppelbergmikosch1997},
Chapter~3.
However, the stationary solution $(Y_i)$ to (\ref{eqkesten}) is not
i.i.d., and therefore one needs to modify (\ref{eqextreme}) as\vadjust{\goodbreak}
follows: the limit
has to be replaced by $\Phi_\alpha^\theta$ for some constant
$\theta\in(0,1)$, the so-called \textit{extremal index} of the sequence
$(Y_i)$; see de Haan et
al.~\cite{dehaanresnickrootzendevries1989};
cf.~\cite{embrechtskluppelbergmikosch1997}, Section 8.4.

The main objective of this paper is to derive
another result which is a consequence
of the power law tails of the marginal distribution of the sequence $(Y_i)$:
we will prove large deviation results for the partial sum sequence
\[
\S_n =Y_1+\cdots+Y_n,\qquad n\ge1,\qquad
\S_0=0.
\]
This means we will derive exact asymptotic results for the
left and right tails of the partial sums $\S_n$. Since we want to
compare these results with those for an i.i.d.
sequence, we recall the corresponding classical results due to A. V. and
S. V. Nagaev~\cite{nagaev1969,nagaev1979} and
Cline and Hsing~\cite{clinehsing1998}.
%
%th1.2 #&#
\begin{theorem}\label{thmnagaev}
Assume that $(Y_i)$ is an i.i.d. sequence with a regularly varying distribution,
that is, there exists an $\alpha>0$, constants $p,q\ge0$ with $p+q=1$
and a slowly varying function $L$ such that
%
%e1.5 #&#
\begin{equation}
\label{eqtailbalance} \pb\{Y>x\} \sim p \frac{L(x)}{x^\alpha}
\quad\mbox{and}\quad \pb\{Y\le-x
\} \sim q \frac{L(x)}{x^\alpha} \qquad\mbox{as $\xto$.}
\end{equation}
\item
Then the following relations
hold for $\alpha>1$ and suitable sequences $b_n\uparrow\infty$:
%
%e1.6 #&#
\begin{equation}
\label{eqlda} \lim_{\nto} \sup_{x\ge b_n} \biggl\llvert
\frac{\pb\{\S_n- \eb\S_n>x\}}{n \pb\{|Y|>x\}}- p\biggr\rrvert =0
\end{equation}
and
%
%e1.7 #&#
\begin{equation}
\label{eqldb} \lim_{\nto} \sup_{x\ge b_n}\biggl\llvert
\frac{\pb\{\S_n- \eb\S_n\le-x\}}{n \pb\{
|Y|>x\}}- q\biggr\rrvert =0.
\end{equation}
If $\alpha>2$ one can choose $b_n=\sqrt{a n \log n}$, where
$a>\alpha-2$, and for $\alpha\in(1,2]$, $b_n=n^{\delta+1/\alpha}$
for any $\delta>0$.
\item
For $\alpha\in(0,1]$, (\ref{eqlda}) and (\ref{eqldb}) remain
valid if
the centering
$\eb\S_n$ is replaced by~$0$ and $b_n=n^{\delta+1/\alpha}$
for any $\delta>0$.
\end{theorem}
For $\alpha\in(0, 2]$ one can choose a smaller bound $b_n$
if one knows the
slowly varying function $L$ appearing in (\ref{eqtailbalance}).
%Moreover,
%if $\alpha\in(0,1)$ or $\alpha=1$ and $\eb|Y|=\infty$,
%a result similar to equations
%constants $\eb\S_n$ replaced by zero.
A functional version of
Theorem~\ref{thmnagaev} with multivariate regularly varying
summands was proved in Hult
et al.~\cite{hultlindskogmikoschsamorodnitsky2005} and the results
were used to prove asymptotic results about multivariate ruin probabilities.
Large deviation results for i.i.d. heavy-tailed summands are
also known when the distribution of the summands is subexponential, including
the case of regularly varying tails; see the
recent paper by Denisov et al.~\cite{denisovdiekershneer2008} and the
references therein. In this case, the regions where the large
deviations hold
very much depend on the decay rate of the tails of the summands. For
semi-exponential tails (such as for the log-normal and the
heavy-tailed Weibull distributions) the large deviation regions
$(b_n,\infty)$ are much
smaller than those for summands with regularly varying tails. In particular,
$x=n$ is not necessarily contained in $(b_n,\infty)$.\vadjust{\goodbreak}

The aim of this paper is to study large deviation probabilities for a particular
dependent sequence $(Y_n)$ as described in Kesten's Theorem~\ref{thmkesten}.
For dependent sequences $(Y_n)$ much less is known
about the large deviation probabilities for the partial sum process $(\S_n)$.
Gantert~\cite{gantert2000} proved large deviation results of
logarithmic type
for mixing subexponential random variables.
Davis and Hsing~\cite{davishsing1995} and
Jakubowski~\cite{jakubowski1993,jakubowski1997} proved large
deviation results
of the following type: there exist sequences $s_n\to\infty$ such that
\[
\frac{\pb\{\S_n>a_n s_n\}}{n \pb\{Y> a_n s_n\}}\to c_\alpha
\]
for suitable positive constants $c_\a$ under the assumptions that
$Y$ is regularly varying with index $\alpha\in(0,2)$, $n P(|Y|>a_n)\to1$,
and $(Y_n)$ satisfies some mixing conditions. Both Davis and Hsing
\cite{davishsing1995} and
Jakubowski~\cite{jakubowski1993,jakubowski1997}
could not specify the rate at which the sequence $(s_n)$
grows to infinity, and an extension to $\a>2$ was
not possible. These facts limit the
applicability of these results, for example, for deriving the asymptotics
of ruin probabilities for the random walk $(\S_n)$. Large deviations
results for particular stationary sequences $(Y_n)$ with regularly varying
finite-dimensional distributions were proved in Mikosch and Samorodnitsky
\cite{mikoschsamorodnitsky2000} in the case of linear processes
with i.i.d. regularly varying noise and in Konstantinides and Mikosch
\cite{konstantinidesmikosch2005}
for solutions $(Y_n)$ to the stochastic recurrence equation~(\ref
{eqkesten}), where $B$ is
regularly varying with index $\a>1$ and $\eb A^\a<1$. This means that Kesten's
condition (\ref{eqkestena}) is not satisfied in this case, and the
regular variation of $(Y_n)$ is due to the regular variation of $B$.
For these processes,
large deviation results and ruin bounds are easier to derive by
applying the
``heavy-tail large deviation heuristics'': a large value of $\S_n$
happens in the
most likely way, namely it is due to one very large value in the
underlying regularly varying noise
sequence, and the particular dependence structure of the sequence
$(Y_n)$ determines
the clustering behavior of the large values of $\S_n$.
This intuition fails when one deals with the partial sums $\S_n$
under the conditions of Kesten's Theorem~\ref{thmkesten}: here a
large value of $\S_n$ is not due to a single large value of the
$B_n$'s or $A_n$'s but to large values of the products $A_1\cdots A_n$.

The paper is organized as follows. In Section~\ref{sec2} we prove an
analog to Theorem~\ref{thmnagaev} for the partial sum sequence $(\S_n)$
constructed from the solution to the stochastic recurrence equation
(\ref{eqkesten}) under the
conditions of Kesten's Theorem~\ref{thmkesten}. The proof of this
result is rather technical: it is given in Section~\ref{sec3} where
we split the proof into a series of auxiliary results. There we
treat the different cases $\a\le1$, $\a\in(1,2]$ and $\a>2$ by
different tools and methods. In particular, we will use exponential tail
inequalities which are suited for the three distinct situations.
In contrast to the i.i.d. situation described in Theorem~\ref{thmnagaev},
we will show that the $x$-region where the large deviations hold cannot
be chosen
as an infinite interval $(b_n,\infty)$ for a suitable lower bound
$b_n\to\infty$, but one also needs upper bounds $c_n\ge b_n$.
In Section~\ref{ruin} we apply the large deviation results to get
precise asymptotic
bounds for
the ruin probability related to the random walk $(\S_n)$. This ruin bound
is an analog of the celebrated result by Embrechts and Veraverbeke
\cite{embrechtsveraverbeke1982} in the case of a random walk with
i.i.d. step sizes.

%s2 #&#
\section{Main result}\label{sec2}
The following is the main result of this paper.
It is an analog of the well-known large deviation result of Theorem
\ref{thmnagaev}.
%
%th2.1 #&#
\begin{theorem}\label{mthm} Assume that the conditions of
Theorem~\ref{thmkesten} are satisfied and additionally there exists
$\eps>0$
such that $\eb A^{\a+\eps}$ and $\eb|B|^{\a+\eps}$ are finite.
Then the following relations hold:
\begin{longlist}[(1)]
\item[(1)] For $\a\in(0,2] $, $M>2$,
%
%e2.1 #&#
\begin{equation}
\label{unres} \sup_n \sup_{n^{1/\a}(\log n)^M\leq x
}\frac{\pb
\{ \S_n-d_n>x\}}{n \pb\{ |Y|>x\} } <\infty.
\end{equation}
If
additionally $e^{s_n} \ge n^{1/\a}(\log n)^M$ and $\lim_{n\to
\infty}s_n/n=0$, then
%
%e2.2 #&#
\begin{equation}
\label{res} \lim_{\nto}\sup_{n^{1/\a}(\log n)^M\leq x \le
e^{s_n}}\biggl\llvert
\frac{\pb\{ \S_n-d_n>x\}}{n \pb\{ |Y|>x\} } -\frac{c_\infty^+ c_{\infty}}{c_\infty^++c_\infty^-}\biggr\rrvert =0,
\end{equation}
where $d_n=0$ or $d_n=\eb\S_n$ according as $\alpha\in
(0,1]$ or $\alpha\in(1,2]$.
\item[(2)]
For $\a>2$ and any $c_n\to\infty$,
%
%e2.3 #&#
\begin{equation}
\label{2unres} \sup_n\sup_{ c_n n^{0.5}\log n\leq x }\frac{\pb\{ \S_n-\eb\S_n>x\}}{n \pb\{ |Y|>x\} } <\infty.
\end{equation}
If
additionally $c_nn^{0.5}\log n\le e^{s_n}$ and $\lim_{n\to
\infty}s_n/n=0$, then
%
%e2.4 #&#
\begin{equation}
\label{2res} \lim_{\nto}\sup_{c_n
n^{0.5}\log n\leq x \leq e^{s_n}}\biggl\llvert
\frac{\pb\{ \S_n-\eb\S_n>x\}}{n \pb\{ |Y|>x\} } -\frac{c_\infty^+
c_{\infty}}{c_\infty^++c_\infty^-}\biggr\rrvert =0.
\end{equation}
\end{longlist}
\end{theorem}
Clearly, if we exchange the variables $B_n$ by $-B_n$ in the
above results we obtain the corresponding asymptotics for the left tail
of $\S_n$. For example, for $\alpha>1$ the following relation
holds uniformly for the $x$-regions indicated above:
\[
\lim_{n\to\infty}\frac{\pb\{ \S_n-n\eb Y \le-x\} }{n \pb\{ |Y|>x\} } =\frac{c_\infty^-c_{\infty}}{c_\infty^++c_\infty^-}.
\]

%re2.2 #&#
\begin{remark}
The deviations of Theorem~\ref{mthm} from the i.i.d. case (see
Theorem~\ref{thmnagaev})
are two-fold. First, the extremal clustering in the sequence $(Y_n)$
manifests in the presence of the additional constants $c_\infty$
and $c_\infty^\pm$. Second, the precise large deviation bounds (\ref{res})
and (\ref{2res}) are proved for $x$-regions bounded from above
by a sequence $e^{s_n}$ for some $s_n\to\infty$ with \mbox{$s_n/n\to0$}.
%It is shown in the course of the proof (see Section~\ref{counter})
%that \eqref{res} and \eqref{2res} cannot be extended to unbounded
%$x$-regions, i.e. in the latter case the upper bounds \eqref{unres} and
Mikosch and Wintenberger\vadjust{\goodbreak}
\cite{mikoschwintenberger2011} extended Theorem~\ref{mthm} to more
general classes of stationary sequences $(Y_t)$. In particular, they proved
similar results for stationary Markov chains with regularly varying
finite-dimensional distributions,
satisfying a
drift condition. The solution $(Y_t)$ to (\ref{eqkesten})
is a special case of this setting if the distributions of $A,B$ satisfy
some additional conditions. Mikosch and Wintenberger \cite
{mikoschwintenberger2011} use a regeneration argument to explain that
the large deviation results
do not hold uniformly in the unbounded
$x$-regions $(b_n,\infty)$ for suitable sequences $(b_n)$, $b_n\to\infty$.
\end{remark}

%s3 #&#
\section{Proof of the main result}\label{sec3}
%s3.1 #&#
\subsection{Basic decompositions}\label{subsecbasic}
In what follows, it will be convenient to use the following
notation:
\[
\Pi_{ij}= \cases{A_i\cdots A_j, &\quad $i\le
j$,
\cr
1, &\quad otherwise,} \quad\mbox{and}\quad \Pi_j=\Pi_{1j}
\]
and
\[
\wt Y_i=\Pi_{2i} B_1+
\Pi_{3i}B_2+\cdots+\Pi_{ii}B_{i-1}+B_i,\qquad i\ge1.
\]
Since
$Y_i=\Pi_i Y_0+\wt Y_i$,
the following decomposition is straightforward:
%
%e3.1 #&#
\begin{equation}
\label{eqg1} \S_n = Y_0 \sum
_{i=1}^n\Pi_i+\sum
_{i=1}^n\wt Y_i=: Y_0
\eta_n+\wt\S_n,
\end{equation}
where
%
%e3.2 #&#
\begin{equation}
\label{tsn} \wt\S_n = \wt Y_1 + \cdots+ \wt
Y_n \quad\mbox{and}\quad \eta_n=\Pi_1+\cdots+
\Pi_n,\qquad n\ge1.
\end{equation}
%
%We start with some rough bound on the tail of $|Y_0|\eta_n$.
In view of (\ref{eqg1}) and Lemma~\ref{lemnew} below
it suffices to bound the ratios
\[
\frac{\pb\{ \wt
\S_n-\wt d_n>x\}}{n \pb\{ |Y|>x\} }
\]
uniformly for the considered $x$-regions, where $\wt d_n = \eb\wt\S_n$
for $\alpha>1$ and $\wt d_n = 0$ for $\alpha\le1$.

The proof of the following bound is given at the end of this subsection.
%
%le3.1 #&#
\begin{lemma}\label{lemnew}Let $(s_n)$ be a sequence such that
$s_n/n\to
0$. Then
for any sequence $(b_n)$ with $b_n\to\infty$ the following relations
hold:
\[
\lim_{\nto} \sup_{ b_n\le x\le
e^{s_n}} \frac{\pb\{ |Y_0| \eta_n>x\}}{n \pb\{|Y|>x\}}=0
\quad\mbox{and}\quad
\limsup_{\nto}\sup_{b_n\le x} \frac{\pb\{ |Y_0| \eta_n>x\}}{n
\pb\{|Y|>x\}}<\infty.
\]
\end{lemma}
Before we further decompose $\wt\S_n$ we introduce some notation to
be used throughout the proof.
For any $x$ in the considered large deviation regions:
\begin{itemize}
\item$m=[(\log x)^{0.5 +\sigma}]$ for some positive number
$\sigma<1/4$, where \mbox{$[\cdot]$} denotes the integer part.
\item$n_0=[\rho^{-1}\log x]$, where
$\rho= \eb(A^{\a}\log A)$.
\item$n_1=n_0-m$ and $n_2=n_0+m$.\vadjust{\goodbreak}
\item For $\a>1$, let $D$ be the smallest integer such that $-D\log
\eb A >\a-1$. Notice that the latter inequality makes sense since
$\eb A<1$ due to (\ref{eqkestena})
and the convexity
of the function $\psi(h)=\eb A^{h}$, $h>0$.
\item For $\a\leq1$, fix some $\b<\a$, and let $D$ be the
smallest integer such that $-D \log\eb A^{\b}>\a-\b$ where, by the
same remark as above, $\eb A^{\b}<1$.
\item Let $n_3$ be the smallest integer satisfying
%
%e3.3 #&#
\begin{equation}
\label{eq33} D\log x\le n_3,\qquad x>1.
\end{equation}
Notice that since the function $\Psi(h) = \log\psi(h)$ is convex,
putting $\b=1$ if \mbox{$\a>1$}, by the choice of $D$ we have $\frac
1D<\frac{\Psi(\a)-\Psi(\b)}{\a-\b} < \Psi'(\a)=\rho$;
therefore $n_2<n_3$ if $x$ is sufficiently large.
\end{itemize}
For fixed $n$, we change the indices
$i \to j=n-i+1$ and, abusing notation and suppressing the dependence
on $n$, we reuse the notation
\[
\wt Y_j= B_j+\Pi_{jj}B_{j+1}+
\cdots+\Pi_{j,n-1}B_n.
\]
Writing
$n_4=\min(j+n_3,n)$, we further decompose $\wt Y_j$,
%
%e3.4 #&#
\begin{equation}
\label{eq2} \wt Y_j=\wt U_j+\wt W_j=
B_j+\Pi_{jj}B_{j+1}+\cdots+\Pi_{j,n_4-1}
B_{n_4}+\wt W_j.
\end{equation}
Clearly, $\wt W_j$ vanishes if $j\geq n-n_3$
and therefore the following lemma is nontrivial only for
$n>n_3$. The proof is given at the end of this subsection.
%
%le3.2 #&#
\begin{lemma}\label{red2} For any
small $\delta>0$, there exists a constant $c>0$ such that
%
%e3.5 #&#
\begin{equation}
\label{eq44} \pb \Biggl\{ \Biggl|\sum_{j=1}^{n}(
\wt W_j- c_{j}) \Biggr|>x \Biggr\}\le c n x^{-\a-\delta},\qquad x>1,
\end{equation}
where $c_j=0$ or $c_j=\eb\wt W_j$ according as $\alpha\le1$ or
$\alpha>1$.
\end{lemma}
By virtue of (\ref{eq44}) and (\ref{eq2}) it
suffices to study the probabilities\break $ \pb \{ \sum_{j=1}^{n}(\wt
U_j-a_j)>x \}$, where $a_j=0$ for $\alpha\le1$ and $a_j=\eb\wt
U_j$ for $\alpha>1$.

We further decompose $\wt U_i$ into
%
%e3.6 #&#
\begin{equation}
\label{dec2} \wt U_i= \wt X_i + \wt
S_i +\wt Z_i,
\end{equation}
where for $i\le n-n_3$,
%
%e3.7 #&#
\begin{eqnarray}
\label{eq4a} \wt X_i&=&B_i+\Pi_{ii}B_{i+1}+
\cdots+\Pi_{i,i+n_1-2}B_{i+n_1-1},
\nonumber
\\
\wt S_i&=&\Pi_{i,i+n_1-1}B_{i+n_1}+\cdots +
\Pi_{i,i+n_2-1}B_{i+n_2},%\label{eq6}
\\
\wt Z_i&=&\Pi_{i,i+n_2}B_{i+n_2+1}+\cdots +
\Pi_{i,i+n_3-1}B_{i+n_3}.
\nonumber
\end{eqnarray}
For $i>n-n_3$, define $\wt X_i,
\wt S_i, \wt Z_i $ as follows: for $n_2<n-i<n_3$ choose
$\wt X_i, \wt S_i$ as above and
\[
\wt Z_i =\Pi_{i,i+n_2}B_{i+n_2+1}+\cdots+
\Pi_{i,n-1}B_{n}.\vadjust{\goodbreak}
\]
For $n_1\le n-i\le n_2$, choose $\wt Z_i=0$, $\wt X_i$ as before
and
\[
\wt S_i=\Pi_{i,i+n_1-1}B_{i+n_1}+\cdots+
\Pi_{i,n-1}B_{n}.
\]
Finally, for $n-i<n_1$, define $\wt S_i=0, \wt Z_i=0$ and
\[
\wt X_i=B_i+\Pi_{ii}B_{i+1}+
\cdots+\Pi_{i,n-1}B_{n}.
\]

Let $p_1,p,p_3$ be the
largest integers such that $p_1n_1\leq n-n_1+1$, $pn_1\leq n-n_2$ and
$p_3n_1\leq n-n_3$, respectively. We study the asymptotic tail behavior of
the corresponding block sums given by
%
%e3.8 #&#
\begin{equation}
\label{eq4}\qquad X_j=\sum_{i=(j-1)n_1+1}^{jn_1}
\wt X_i,\qquad S_j=\sum_{i=(j-1)n_1+1}^{jn_1}
\wt S_i,\qquad Z_j=\sum_{i=(j-1)n_1+1}^{jn_1}
\wt Z_i,
\end{equation}
where $j$ is less or equal
$p_1, p, p_3$, respectively.

The remaining steps of the proof are organized as follows:
\begin{itemize}
\item
Section~\ref{subsecxz}. We show that the $X_j$'s and $Z_j$'s
do not contribute to the considered large deviation probabilities.
This is the content of Lemmas~\ref{x} and~\ref{z}.
\item
Section~\ref{subsecs}. We provide bounds for the
tail probabilities of $S_j$; see Proposition~\ref{block} and Lemma~\ref{prod}.
These bounds are the main ingredients in the proof of the
large deviation result.
\item
Section~\ref{subsecfinal}. In Proposition~\ref{pomoc} we combine the
bounds provided in the previous subsections.
\item
Section~\ref{subsecfinala}: we apply Proposition~\ref{pomoc} to
prove the main result.
\end{itemize}
\begin{pf*}{Proof of Lemma~\ref{lemnew}}
The infinite series $\eta=\sum_{i=0}^\infty\Pi_i$ has the
distribution of the stationary solution to the stochastic recurrence
equation (\ref{eqkesten})
with \mbox{$B\equiv1$} a.s., and therefore, by Theorem~\ref{thmkesten},
$
P(\eta>x)\sim c_\infty x^{-\alpha}, \xto.
$
It follows from a
slight modification of Jessen and Mikosch
\cite{jessenmikosch2006}, Lemma~4.1(4), and the independence of
$Y_0$ and $\eta$ that
%
%e3.9 #&#
\begin{equation}
\label{jm} \pb\bigl\{|Y_0| \eta>x\bigr\}\sim c x^{-\alpha
}\log x,\qquad \xto.
\end{equation}
Since $s_n/n\to0$ as $\nto$ we have
\[
\sup_{
b_n\le x\le e^{s_n}} \frac{\pb\{ |Y_0| \eta_n>x\}}{n \pb\{|Y|>x\}
}\le\sup_{ b_n\le x\le
e^{s_n}}
\frac{\pb\{ |Y_0| \eta>x\}}{n \pb\{|Y|>x\}}\to0.
\]
There exist $c_0,x_0>0$ such that
$P\{|Y_0|>y\}\le c_0y^{-\alpha}$ for $y>x_0$. Therefore
\[
\pb\bigl\{ |Y_0| \eta_n>x\bigr\} \le\P\{x/\eta_n\le
x_0\}+ c_0x^{-\a}\eb\eta_n^\a
\one_{\{x/\eta_n>x_0\}}\le c x^{-\a}\eb\eta_n^\a.
\]
By Bartkiewicz et al.
\cite{bartkiewiczjakubowskimikoschwintenberger2009},
$\eb\eta_n^{\a}\le cn$. Hence
\[
I_n=\sup_{b_n\le x} \frac{\pb\{ |Y_0| \eta_n>x\}}{n \pb\{|Y|>x\}
}\le\sup_{ b_n\le x }
\frac{cx^{-\a}\eb\eta_n^{\a}}{n \pb\{|Y|>x\}}<\infty.
\]
This concludes the proof.\vadjust{\goodbreak}
\end{pf*}
\begin{pf*}{Proof of Lemma~\ref{red2}}
Assume first that $\a>1$. Since $\eb\wt W_j$ is
finite, $-D\log
\eb A >\a-1$ and $D\log x\le n_3$, we have for some positive
$\d$
%
%e3.10 #&#
\begin{equation}
\label{eqr1} \eb|\wt W_j|\leq\frac{(\eb A)^{n_3}}{1-\eb A}\eb|B| \leq c
e^{D \log x \log\eb A} \leq c x^{-(\a-1)-\d},
\end{equation}
and hence by Markov's inequality
\[
\pb \Biggl\{ \Biggl| \sum_{j=1}^{n}(\wt
W_j-\eb\wt W_j) \Biggr| >x \Biggr\}\leq 2 x^{-1}
\sum_{j=1}^{n}\eb|\wt W_j|
\leq c n x^{-\a-\d}.
\]
If $\beta<\a\le1$ an application of Markov's inequality yields
for some positive~$\d$,
\begin{eqnarray*}
\pb \Biggl\{ \sum_{j=1}^{n}\wt
W_j>x \Biggr\} &\leq& x^{-\beta}\sum
_{j=1}^{n}\eb|\wt W_j|^{\b}
\leq x^{-\beta}\frac{n\eb
|B|^{\b}(\eb A^{\b})^{n_3}}{(1-\eb
A^{\b})}
\\
&\le& c x^{-\beta} n e^{D \log x \log\eb A^{\b}} \leq c n x^{-\a-\d}.
\end{eqnarray*}
In the last step we used the fact that $-D \log\eb A^{\b}>\a-\b$.
This concludes the proof of the lemma.
\end{pf*}
%
%s3.2 #&#
\subsection{\texorpdfstring{Bounds for $\P\{X_j>x\}$ and $\P\{Z_j>x\}$}
{Bounds for $P\{X_j>x\}$ and $P\{Z_j>x\}$}}\label{subsecxz}
%In a series of auxiliary results
We will now study the tail behavior of
the single block sums
$X_1,Z_1$ defined in (\ref{eq4}). We start with a useful auxiliary result.
%
%le3.3 #&#
\begin{lemma}\label{exp} Assume
$\psi(\alpha+\epsilon)=\eb A^{\a+\epsilon}<\infty$
for some $\epsilon>0$. Then there is a constant $C=C(\epsilon)>0$
such that
$
\psi(\a+ \gamma)\leq C e^{\rho\gamma}$
for $|\gamma|\le\epsilon/2$,
where $\rho=\eb(A^\alpha\log A)$.
\end{lemma}
\begin{pf}
By a Taylor expansion and since $\psi(\alpha)=1$, $\psi'(\a)=\rho$,
we have for some $\theta\in(0,1)$,
%
%e3.11 #&#
\begin{equation}
\label{eqtaylor} \psi(\a+\gamma)=1+\rho\gamma+ 0.5 \psi''(
\a+\theta\gamma) \gamma^2.
\end{equation}
If $|\theta\gamma| <\epsilon/2$, then, by assumption,
$\psi''(\a+\theta\gamma)=\eb A^{\a+\theta\gamma}(\log A)^2$ is
bounded by a constant $c>0$. Therefore,
\[
\psi(\a+\gamma)\leq1+\rho\gamma+c\gamma^2 =e^{\log(1+\rho\gamma+c \gamma^2)}\leq C
e^{\rho\gamma}.
\]
\upqed\end{pf}
%
%We notice that $|\wt X_1|$ and $|\wt Z_1|$ are stochastically
%dominated by $\underline{ X}_1$ and $\underline{Z}_1$, respectively.
%Therefore the bounds in Lemmas~\ref{x} and~\ref{z} also apply to the
%tails of $|\wt X_1|$ and $|\wt Z_1|$, respectively.
The following lemma ensures that the $X_i$'s do not contribute to the
considered large deviation probabilities.
%
%le3.4 #&#
\begin{lemma}\label{x} There exist positive constants
$C_1,C_2,C_3 $ such that
\[
\pb\{X_1 > x\} \le\pb\{ \underline{X}_1 >x\} \leq
C_1 x^{-\a
}e^{-C_2 (\log
x)^{C_3}},\qquad x>1,
\]
where
\[
\underline{ X}_1=\sum_{i=1}^{n_1}\bigl(|B_i|+
\Pi_{ii}|B_{i+1}|+\cdots +\Pi_{i,i+n_1-2}
|B_{i+n_1-1}|\bigr).
\]
\end{lemma}
\begin{pf} We have
$\underline{X}_1= \sum_{k=m+1}^{n_0}R_k$,
where for $m< k\leq n_0$,
\begin{eqnarray*}
R_k&=&\Pi_{1,n_0-k}|B_{n_0-k+1}|+\cdots +
\Pi_{i,i+n_0-k-1}|B_{i+n_0-k}|+\cdots \\
&&{}+\Pi_{n_1,n_1+n_0-k-1}|B_{n_1+n_0-k}|.
\end{eqnarray*}
Notice that for $x$ sufficiently large,
\[
\Biggl\{\sum_{k=m+1}^{n_0}R_k >x
\Biggr\} \subset\bigcup_{k=m+1}^{n_0}\bigl\{
R_k > x/k^3\bigr\}.
\]
Indeed, on the set $\{R_k\leq
x/k^3, m<k\le n_0\}$ we have for some $c>0$ and sufficiently
large $x$, by the definition of $m=[(\log x)^{0.5 +\sigma}]$,
\[
\sum_{k=m+1}^{n_0}R_k\leq
\frac{x}{m+1}\sum_{k=1}^{\infty}
\frac{1} {
k^2}\leq c \frac{x}{(\log x)^{0.5 +\sigma}}<x.
\]
We conclude
that, with $I_k=\pb\{R_k > x/k^3\}$,
\[
\pb \Biggl\{\sum_{k=m+1}^{n_0}R_k
>x \Biggr\} \leq\sum_{k=m+1}^{n_0}I_k.
\]
Next we study the probabilities $I_k$. Let $\d= (\log x)^{-0.5}$. By
Markov's inequality,
\[
I_k \leq\bigl(x/k^3\bigr)^{-(\a+\d)} \eb
R_k^{\a+\d} \leq\bigl(x/k^3\bigr)^{-(\a+\d)}
n_0^{\a+\d}\bigl(\eb A^{\a+\d}\bigr)^{n_0-k}
\eb |B|^{\a+\d}.
\]
By Lemma~\ref{exp} and the definition of $n_0=[\rho^{-1}\log x]$,
\[
I_k \leq c \bigl(x/k^3\bigr)^{-(\a+\d)}
n_0^{\a+\d
} e^{(n_0-k)\rho\d} \leq c x^{-\a}
k^{3(\a+\d)} n_0^{\a+\d} e^{-k\rho\d}.
\]
Since $k\geq(\log x)^{0.5+\sigma}\ge m$ there are positive constants
$\zeta_1,\zeta_2$ such that
$k\d\geq k^{\zeta_1}(\log x)^{\zeta_2}$ and therefore for
sufficiently large $x$ and appropriate positive constants $C_1,C_2,C_3$,
\[
\sum_{k=m+1}^{n_0}I_k \leq c
x^{-\a} n_0^{\a+\d} \sum
_{k=m+1}^{n_1} e^{-\rho k^{\zeta_1} (\log x)^{\zeta_2} }k^{3(\a
+\d)} \leq
C_1 x^{-\a} e^{-C_2 (\log x)^{C_3}}.
\]
This finishes the proof.
\end{pf}
The following lemma ensures that the $Z_i$'s do not contribute to the
considered large deviation probabilities.
%
%le3.5 #&#
\begin{lemma}\label{z} There exist positive constants $C_4,C_5,C_6$
such that
\[
\pb\{ Z_1 >x\} \le \pb\{ \underline{Z}_1 >x\} \leq
C_4 x^{-\a}e^{-C_5 (\log
x)^{C_6}},\qquad x> 1,
\]
where
\[
\underline{Z}_1=\sum_{i=1}^{n_1}\bigl(
\Pi_{i,i+n_2} |B_{i+n_2+1}| +\cdots+\Pi_{i,i+n_3-1}|B_{i+n_3}|\bigr).
\]
\end{lemma}
\begin{pf} We have
$\underline{Z}_1 = \sum_{k=1}^{n_3-n_2}\wt R_k$, where
\begin{eqnarray*}
\wt R_k&=&\Pi_{1,n_2+k}|B_{n_2+k+1}|+\cdots+
\Pi_{i,i+n_2+k-1}|B_{i+n_2+k}|+\cdots\\
&&{}+ \Pi_{n_1,n_1+n_2+k-1}|B_{n_1+n_2+k}|.
\end{eqnarray*}
As in the proof of
Lemma~\ref{x} we notice that, with
$J_k=\pb
\{ \wt R_k> x/(n_2+k)^3\}$, for $x$ sufficiently large,
\[
\pb\Biggl\{ \sum_{k=1}^{n_3-n_2}\wt
R_k >x\Biggr\} \leq\sum_{k=1}^{n_3-n_2}J_k.
\]
Next we study the probabilities $J_k$.
Choose $\d=(n_2+k)^{-0.5}<\epsilon/2$
with $\epsilon$ as in Lemma~\ref{exp}. By Markov's inequality,
\[
J_k \leq\bigl((n_2+k)^{3}/x
\bigr)^{\a-\d} \eb\wt R_k^{\a-\d} \leq
\bigl((n_2+k)^{3}/x\bigr)^{\a-\d}
n_1^{\a-\d} \bigl(\eb A^{\a
-\d}\bigr)^{n_2+k}
\eb|B|^{\a-\d}.
\]
By Lemma~\ref{exp} and since $n_2+k=n_0+m+k$,
\[
\bigl(\eb A^{\a
-\d}\bigr)^{n_2+k}\leq c e^{-\d\rho(n_2+k)} \leq c
x^{-\d} e^{-\d\rho(m+k)}.
\]
There is $\zeta_3>0$ such that
$
\d(m+k)\geq(\log x +k)^{\zeta_3}$. Hence, for appropriate constants
$C_4,C_5, C_6 >0$,
\[
\sum_{k=1}^{n_3-n_2}J_k\leq c
x^{-\alpha} n_1^{\a-\d}\sum
_{k=1}^{n_3-n_2}(n_2+k)^{3(\a-\d
)}e^{-\rho
(\log x+k)^{\zeta_3}}
\leq C_4 x^{-\a}e^{-C_5 (\log x)^{C_6}}.
\]
This finishes the proof.
\end{pf}
%
%s3.3 #&#
\subsection{\texorpdfstring{Bounds for $\P\{S_j>x\}$}{Bounds for $P\{S_j>x\}$}}\label{subsecs}
The next proposition is a first major step toward the proof of the main
result. For the formulation of the result and its proof, recall the
definitions of $\wt S_i$ and $S_i$ from
(\ref{eq4a}) and (\ref{eq4}), respectively.

%pr3.6 #&#
\begin{proposition}\label{block} Assume that $c_\infty^+>0$ and let
$(b_n)$ be
any sequence such that $b_n\to\infty$. Then the following relation holds:
%
%e3.12 #&#
\begin{equation}
\label{es} \lim_{n\to\infty}\sup_{x\ge b_n} \biggl|\frac
{\pb\{
S_1>x\} }{n_1\pb\{ Y>x\}}
-c_{\infty} \biggr|=0.
\end{equation}
If
$c_\infty^+=0$, then
%
%e3.13 #&#
\begin{equation}
\label{ess} \lim_{n\to\infty}\sup_{x\ge
b_n}\frac{\pb\{ S_1>x\} }{n_1\pb\{ |Y|>x\}}=0.
\end{equation}
\end{proposition}
The proof depends on the following auxiliary result whose proof is
given in Appendix~\ref{app2}.
%
%le3.7 #&#
\begin{lemma}\label{estim} Assume that $Y$ and
$\eta_k$ [defined in (\ref{tsn})] are independent and
$\psi(\a+\epsilon)=\eb A^{\a+ \epsilon}<\infty$ for some
$\epsilon>0$. Then for
$n_1=n_0-m=[\rho^{-1}\log x]- [(\log x)^{0.5+\sigma}]$ for some
$\sigma<1/4$ and\vadjust{\goodbreak}
any sequences $b_n\to\infty$ and $r_n\to\infty$
the following relation holds:
\[
\lim_{n\to
\infty}\sup_{r_n \le k\le n_1, b_n\le x} \biggl|\frac{\pb\{ \eta_k Y>x\}}{k \pb\{ Y>x\}}-c_{\infty}
\biggr|=0,
\]
provided $c_\infty^+>0$. If
$c_\infty^+=0$, then
\[
\lim_{n\to\infty}\sup_{r_n\le k\le n_1,
b_n\le x} \frac{\pb\{ \eta_k Y>x\}}{k \pb\{ |Y|>x\}}=0.
\]
\end{lemma}
%
%and \eqref{eqt01}
%follow
%by observing that $U\eqd\eta_{n_1}$ and $Y\eqd T_1+T_2$. \ere
%
\begin{pf*}{Proof of Proposition~\ref{block}}
For $i\leq n_1$, consider
\begin{eqnarray*}
\wt S_i+S'_i
&=& \Pi_{i,n_1}B_{n_1+1}+\cdots+\Pi_{i,i+n_1-2}
B_{i+n_1-1}+\wt S_i +\Pi_{i,i+n_2}B_{i+n_2+1}+
\cdots\\
&&{}+\Pi_{i,n_2+n_1-1}B_{n_2+n_1}
\\
&=&\Pi_{i,n_1}(B_{n_1+1}+A_{n_1+1}B_{n_1+2}+
\cdots+\Pi_{n_1+1,n_2+n_1-1}B_{n_2+n_1}).
\end{eqnarray*}
Notice that
\[
\pb\bigl\{ \bigl|S_1'+\cdots+S'_{n_1}\bigr|>x
\bigr\} \leq n_1 \pb\bigl\{ \bigl|S_1'\bigr|>x/n_1
\bigr\}.
\]
Therefore and by virtue of Lemmas~\ref{x} and~\ref{z},
%(cf. Remark~\ref{rem2})
there exist
positive constants $C_7,C_8,C_9$ such that
\[
\pb\bigl\{ \bigl|S_1'+\cdots+S'_{n_1}\bigr|>x
\bigr\} \leq C_7 x^{-\a} e^{-C_8(\log x)^{C_9 }},\qquad x\ge1.
\]
Therefore and since $S_1=\sum_{i=1}^{n_1}\wt S_i$
it suffices for (\ref{es}) to
show that
\[
\lim_{n\to\infty}\sup_{x\ge b_n} \biggl| \frac{\pb\{ S_1+\sum_{i=1}^{n_1}S_i'>x\}}{n_1\pb
\{ Y>x\}} -c_{\infty}
\biggr|=0.
\]
We observe that
\[
S_1+\sum_{i=1}^{n_1}S_i'=:UT_1
\quad\mbox{and}\quad T_1+T_2\eqd Y,
\]
where
\begin{eqnarray*}
U&=&\Pi_{1,n_1}+\Pi_{2,n_1}+\cdots+\Pi_{n_1,n_1},
\\
T_1&=& B_{n_1+1}+\Pi_{n_1+1,n_1+1}B_{n_1+2}+
\cdots+\Pi_{n_1+1,n_2+n_1-1}B_{n_2+n_1},
\\
T_2&=&\Pi_{n_1+1,n_2+n_1}B_{n_2+n_1+1}+\Pi_{n_1+1,n_2+n_1+1}
B_{n_2+n_1+2}+\cdots.
\end{eqnarray*}
Since $U=_d \eta_{n_1}$ and $Y=_d T_1 + T_2$, in view of Lemma \ref
{estim} we obtain
\[
\lim_{n\to\infty}\sup_{x\ge b_n} \biggl|\frac{\pb\{ U(T_1+T_2)>x\}
}{n_1 \pb\{ Y>x\}} -c_{\infty}
\biggr|=0,
\]
provided
$c^+_{\infty}>0$ or
\[
\lim_{n\to\infty}\sup_{x\ge b_n}\frac{\pb\{
U(T_1+T_2)>x\} }{n_1 \pb\{| Y|>x\}}=0,
\]
if $ c_{\infty}^+=0$. Thus to prove the proposition it suffices to
justify the existence of some positive constants $C_{10},C_{11},C_{12}$
such that
%
%e3.14 #&#
\begin{equation}
\label{eqt2} \pb\bigl\{ |UT_2|>x\bigr\} \le C_{10}
x^{-\a
}e^{-C_{11} (\log x)^{C_{12}} },\qquad x>1.
\end{equation}
For this purpose we use the same argument as in the proof of
Lemma~\ref{x}. First we write
\[
\pb\bigl\{ |UT_2|>x\bigr\}\leq\sum_{k=0}^{\infty}
\pb\bigl\{ U \Pi_{n_1+1,n_1+n_2+k} |B_{n_1+n_2+k+1}|>x/(\log x+k)^3
\bigr\}.
\]
Write $\d= (\log x+k)^{-0.5}$. Then by Lemma
\ref{exp}, Markov's inequality and since $n_2=n_0+m$,
\begin{eqnarray*}
&&\pb\bigl\{ U \Pi_{n_1+1,n_1+n_2+k} |B_{n_1+n_2+k+1}|>x/(\log x+k)^3
\bigr\}
\\
&&\qquad\leq(\log x+k)^{3(\a-\d)}x^{-(\a-\d)
} \eb U^{\a-\d}\bigl(\eb
A^{\a-\d}\bigr)^{n_2+k}\eb |B|^{\a-\d}
\\
&&\qquad\leq c (\log x+k)^{3(\a-\d)}x^{-(\a-\d)
} e^{-(n_2+k)\rho\d}
\\
&&\qquad \leq c e^{-(m+k)\rho\d}(\log x+k)^{3(\a-\d)}x^{-\a}.
\end{eqnarray*}
There is $\zeta>0$ such that
$
(m+k)\d\geq(\log x +k)^{\zeta}$ and therefore,
\begin{eqnarray*}
\pb\bigl\{ |UT_2|>x \bigr\} &\leq& c x^{-\a}\sum
_{k=0}^{\infty}e^{-(\log x
+k)^{\zeta}\rho}(\log x +k)^{3(\a-\d)}
\\
&\leq& c x^{-\a}e^{-(\log x)^{\zeta}\rho/2 }.
\end{eqnarray*}
This proves (\ref{eqt2}) and the lemma.
\end{pf*}

Observe that if $|i-j|>2$, then $S_i$ and $S_j$ are independent.
For $|i-j|\leq2$ we have the following bound:
%
%le3.8 #&#
\begin{lemma}\label{prod} The following relation holds for some
constant $c>0$:
\[
\sup_{i\ge1, |i-j|\le2}\pb\bigl\{ |S_i|>x, |S_j|>x \bigr\} \leq c
n_1^{0.5} x^{-\a},\qquad x>1.
\]
\end{lemma}
\begin{pf}
Assume without loss of generality that $i=1$ and $j=2,3$. Then we have
\begin{eqnarray*}
|S_1|&\leq& ( \Pi_{1,n_1}+\cdots+ \Pi_{n_1,n_1} )
\\
&&{}\times \bigl(|B_{n_1+1}|+\Pi_{n_1+1,n_1+1}|B_{n_1+2}|+\cdots+
\Pi_{n_1+1, n_1+n_2-1} |B_{n_2+n_1}| \bigr)\\
&=:&U_1T'_1,
\\
|S_2|&\leq& ( \Pi_{n_1+1,2n_1}+\cdots+\Pi_{2n_1,2n_1} )
\\
&&{}\times \bigl(|B_{2n_1+1}|+\Pi_{2n_1+1,2n_1+1}|B_{2n_1+2}|+\cdots +
\Pi_{2n_1+1,2n_1+n_2-1}|B_{2n_1+n_2}| \bigr)\\[-1pt]
&=:&U_2T'_2,
\\[-1pt]
|S_3|&\leq& ( \Pi_{2n_1+1,3n_1}+\cdots+\Pi_{3n_1,3n_1} )
\\[-1pt]
&&{}\times \bigl(|B_{3n_1+1}|+\Pi_{3n_1+1,3n_1+1}|B_{3n_1+2}|+\cdots +
\Pi_{3n_1+1,3n_1+n_2-1}|B_{3n_1+n_2}| \bigr)\\[-1pt]
&=:&U_3T'_3.
\end{eqnarray*}
We
observe that $U_1\eqd\eta_{n_1}$, $U_i$, $i=1,2,3$, are
independent, $U_i$ is independent of $T_i'$ for each $i$ and the
$T_i'$'s have power law tails with index $\alpha>0$. We conclude
from (\ref{es}) that
\begin{eqnarray*}
\pb\bigl\{ |S_1|>x, |S_2|>x\bigr\} &\leq& \pb\bigl\{
T'_1>x n_1^{-1/(2\a)}\bigr\}\\[-1pt]
&&{} + \pb
\bigl\{ T'_1\leq x n_1^{-1/(2\a)},
U_1T'_1>x, U_2T'_2>x
\bigr\}
\\[-1pt]
& \leq&c n_1^{0.5} x^{-\a}+ \pb\bigl\{
n_1^{-1/(2\a
)}U_1>1, U_2T'_2>x
\bigr\}
\\[-1pt]
& \leq&c n_1^{0.5} x^{-\a}+ \pb\bigl\{
U_1>n_1^{1/(2\a)}\bigr\} \pb \bigl\{
U_2T'_2>x\bigr\}
\\[-1pt]
& \leq& c n_1^{0.5} x^{-\a}.
\end{eqnarray*}
In the same way we can bound $\pb\{ |S_1|>t, |S_3|>t\}$. We omit
details.
\end{pf}
%
%s3.4 #&#
\subsection{Semi-final steps in the proof of the main theorem}\label
{subsecfinal}
In the following proposition, we combine the various tail bounds
derived in the previous sections. For this reason,
recall the definitions of $X_i$, $S_i$ and $Z_i$ from
(\ref{eq4}) and that $p_1,p,p_3$ are the
largest integers such that $p_1n_1\leq n-n_1+1$, $pn_1\leq n-n_2$ and
$p_3n_1\leq n-n_3$, respectively.
%
%pr3.9 #&#
\begin{proposition}\label{pomoc} Assume the conditions of
Theorem~\ref{mthm}. In particular, consider the following $x$-regions:
\[
\Lambda_n=\cases{ \bigl(n^{1/\a}(\log n)^M,
\infty\bigr), &\quad for $\alpha\in(0,2]$, $M>2$,
\vspace*{2pt}\cr
\bigl(c_n
n^{0.5}\log n,\infty\bigr), &\quad for $\alpha>2$, $c_n\to
\infty$,}
\]
and introduce a sequence $s_n\to\infty$ such that $e^{s_n}\in\Lambda_n$
and $s_n=o(n)$.
Then the following relations hold:
%
%e3.15 #&#
%e3.16 #&#
%e3.17 #&#
%e3.18 #&#
\begin{eqnarray}
\label{sss}\qquad \frac{ c_\infty^+c_{\infty}}{c_\infty^++c_\infty^-}&\ge& \limsup_{n\to\infty}\sup_{x\in\Lambda_n}
\frac{\pb
\{ \sum_{j=1}^p(S_j-c_j)>x\} }{n \pb\{ |Y|>x\} },
\\[-1pt]
\label{ss} 0&=&\lim_{n\to\infty}\sup_{x\in\Lambda_n, \log x\le s_n
} \biggl|\frac{\pb\{ \sum_{j=1}^p(S_j-c_j)>x\} }{n \pb\{ |Y|>x\}
} -
\frac{c_\infty^+c_{\infty}}{c_\infty^++c_\infty^-} \biggr|,
\\[-1pt]
\label{xx} 0&=&\lim_{n\to\infty}\sup_{x\in\Lambda_n} \frac{\pb \{
|\sum_{j=1}^{p_1}(X_j-e_j)|>x \}}{n \pb\{ |Y|>x\}},
\\[-1pt]
\label{zz} 0&=&\lim_{n\to\infty}\sup_{x\in\Lambda_n} \frac{ \pb \{
|\sum_{j=1}^{p_3}(Z_j-z_j)|>x \}}{n \pb\{ |Y|>x\}},
\end{eqnarray}
where $c_j=e_j=z_j=0$ for $\alpha\le1$ and $c_j=\eb S_j$,
$e_j=\eb X_j$, $z_j=\eb Z_j$ for
$\alpha>1$.\vadjust{\goodbreak}
\end{proposition}
\begin{pf}
We split the proof into the different cases corresponding to \mbox{$\alpha
\le1$}, $\alpha\in(1,2]$ and $\a>2$.\vspace*{9pt}

\textit{The case $1<\a\le2$}.

\textit{Step} 1: \textit{Proof of} (\ref{sss}) \textit{and} (\ref{ss}).
Since $M>2$, we can choose
$\xi$ so small that
%
%e3.19 #&#
\begin{equation}
\label{xi} 2+4\xi< M \quad\mbox{and}\quad \xi< 1/(4\a),
\end{equation}
and we write $y=x/(\log n )^{2\xi}$. Consider the following disjoint
partition of $\Omega$:
\begin{eqnarray*}
\Omega_1&=& \bigcap_{j=1}^p
\bigl\{|S_j|\leq y\bigr\},
\\
\Omega_2&=&\bigcup_{1\le i<k\le p}
\bigl\{|S_i|>y, |S_k|>y\bigr\},
\\
\Omega_3&=&\bigcup_{k=1}^{p}
\bigl\{ |S_k|>y, |S_i|\leq y \mbox{ for all } i\neq k\bigr\}.
\end{eqnarray*}
Then for $A=\{ \sum_{j=1}^p(S_j-c_j)>x\}$,
%
%e3.20 #&#
\begin{equation}
\label{eqdecomp} \pb\{A\}=\pb\{ A\cap\Omega_1\} +\pb\{A\cap
\Omega_2 \}+ \pb\{A\cap\Omega_3 \} =:I_1+I_2+I_3.
\end{equation}
Next we treat the terms $I_i$, $i=1,2,3$, separately.

\textit{Step} 1a: \textit{Bounds for $I_2$}. We prove
%
%e3.21 #&#
\begin{equation}
\label{eqi2} \lim_{n\to\infty} \sup_{x\in\Lambda_n } \bigl(x^\a/n
\bigr) I_2=0.
\end{equation}
We have
\[
I_2\le\sum_{1\le i<k\le p}\pb\bigl\{
|S_i|>y, |S_k|>y\bigr\}.
\]
For $k\geq i+3$, $S_k$ and $S_i$ are independent and then, by (\ref{es}),
\[
\pb\bigl\{ |S_i|>y, |S_k|>y\bigr\}=\bigl(\pb\bigl\{
|S_1|>y\bigr\}\bigr)^2 \leq c \bigl(n_1(y)
\bigr)^2 y^{-2\a},
\]
where $n_1(y)$ is defined in the same way as $n_1=n_1(x)$ with $x$
replaced by $y$. Also notice that $n_1(y)\le n_1(x)$.
For $k=i+1$ or $i+2$, we have by Lemma~\ref{prod}
\[
\pb\bigl\{ |S_i|>y, |S_k|>y\bigr\}\leq c \bigl(n_1(y)
\bigr)^{0.5} y^{-\a}.
\]
Summarizing the above estimates and
observing that (\ref{xi}) holds, {we obtain for $x\in\Lambda_n$,}
\begin{eqnarray*}
I_2&\leq& c \bigl[p^2n_1^2y^{-2\a}+pn_1^{0.5}y^{-\a}
\bigr]
\\
&\leq&c n x^{-\a} \bigl[x^{-\a}n (\log n)^{4\xi\a}+(\log
n)^{2\xi\a}n_1^{-0.5} \bigr]
\\
&\leq& c n x^{-\a} \bigl[(\log n)^{(4\xi-M) \a}+(\log
n)^{2\xi
\a
-0.5 } \bigr].
\end{eqnarray*}
This proves (\ref{eqi2}).\vadjust{\goodbreak}

\textit{Step} 1b: \textit{Bounds for $I_1$.} We will prove
%
%e3.22 #&#
\begin{equation}
\label{eqi1} \lim_{\nto}\sup_{x\in\Lambda_n}\bigl(x^\a/n
\bigr) I_1=0.
\end{equation}
For this purpose, we write $S^y_j=S_j\one_{\{ |S_j|\leq y\}}$ and
notice that
$\eb S_j=\eb S^y_j+\eb S_j\one_{\{ |S_j|> y\} }$.
Elementary computations show that
%
%e3.23 #&#
\begin{equation}
\label{eqg3} |S_1|^{\a}\leq n_1^{\max(\a,1)}(2m+1)^{\max(\a,1)}
\eb|B|^{\a}.
\end{equation}
Therefore
by the H\"older and Minkowski inequalities and by (\ref{es}),
\begin{eqnarray*}
|\eb S_j\one_{\{ |S_j|> y\} }|&\leq&\bigl(\eb|S_j|^{\a}
\bigr)^{1/\a
} \bigl(\pb\bigl\{ |S_j| >y\bigr\} \bigr)^{1-1/\a}
\\
&\leq& c (\log x )^{1.5+\sigma}y^{-\a+1} \bigl(n_1(y)
\bigr)^{1-1/\a}
\\
&\leq& c (\log x )^{1.5+\sigma+2\xi(\a-1) }x^{-\a+1} n_1.
\end{eqnarray*}
Let now $\g>1/\a$ and $n^{1/\a}(\log n)^M\le x\le n^\gamma$.
Since $p n_1\le n$ and (\ref{xi}) holds,
%
%e3.24 #&#
\begin{equation}
\label{eq01} p |\eb S_j\one_{\{
|S_j|> y\} }|\leq c (\log x
)^{1.5+\sigma+2\xi(\a-1)
}x^{-\a+1} n=o(x).
\end{equation}
If $x>n^{\g}$, then
\[
x>(\log x )^Mn^{1/\a} \quad\mbox{and}\quad x^{-\a}<(
\log x )^{-M\a}n^{-1 }.
\]
Hence
%
%e3.25 #&#
\begin{equation}
\label{eq02} p |\eb S_j\one_{\{ |S_j|> y\}
}|\leq c x (\log x
)^{1.5+\sigma+2\xi(\a-1) }(\log x )^{-M\a} =o(x).
\end{equation}
Using the bounds (\ref{eq01}) and (\ref{eq02}), we see that for $x$
sufficiently
large,
%
%e3.26 #&#
\begin{eqnarray}\label{eq100}\quad
I_1&\leq& \pb \Biggl\{ \Biggl| \sum_{j=1}^p
\bigl(S^y_j-ES^y_j\bigr)\Biggr|>0.5
x \Biggr\}
\nonumber
\\
&=& \pb \biggl\{ \biggl| \biggl( \sum_{1\le j\le p, j\in\{1,4,7,\ldots\}}+ \sum
_{1\le j\le p, j\in\{2,5,8,\ldots\}}\nonumber\\[-8pt]\\[-8pt]
&&\hspace*{96pt}{}+ \sum_{1\le j\le p, j\in\{3,6,9,\ldots\}}
\biggr) \bigl(S^y_j-ES^y_j\bigr)
\biggr|>0.5 x \biggr\}
\nonumber\\
&\le& 3 \pb \biggl\{ \biggl|\sum_{1\le j\le p, j\in\{1,4,7,\ldots\}
}
\bigl(S^y_j-ES^y_j\bigr) \biggr| >x/6
\biggr\}.\nonumber
\end{eqnarray}
In the last step, for the ease of presentation, we slightly abused
notation since the number of summands in the
3 partial sums differs by a bounded number of terms which, however, do not
contribute to the asymptotic tail behavior of $I_1$. Since the summands
$S_1^y,S_4^y,\ldots$ are i.i.d. and bounded, we may apply Prokhorov's
inequality (\ref{prokh})
to the random variables $R_k=S^y_k-\eb S^y_1$ in
(\ref{eq100})
with $y= x/ (\log n)^{2\xi}$ and $B_p=p \var(S^y_1)$. Then\vadjust{\goodbreak}
$a_n=x/(2y)=0.5 (\log n)^{2\xi}$ and since, in view of (\ref{eqg3}),
$\var(S_1^y)\le y^{2-\a} \eb|S_1|^\alpha$,
\[
I_1\leq c \biggl(\frac{p \var(S^y_1)}{xy} \biggr)^{a_n} \leq c
\bigl((\log n)^{(1.5 +\sigma) \a+2\xi(\a-1)-1} \bigr)^{a_n} \biggl(\frac{n}{x^\a}
\biggr)^{a_n}.
\]
Therefore, for {$x\in\Lambda_n$,}
\[
\bigl(x^\a/n\bigr) I_1 \leq c (\log
n)^{((1.5 +\sigma)\a+2\xi(\a-1))a_n -M\a(a_n-1)},
\]
which tends to zero if $M>2$, $\xi$ satisfies (\ref{xi}) and $\sigma
$ is sufficiently small.

\textit{Step} 1c: \textit{Bounds for $I_3$.} Here we assume $c_\infty^+>0$. In
this case,
we can bound $I_3$ by using the following inequalities:
for every $\epsilon
>0$, there is $n_0$ such that for $n\ge n_0$, uniformly for {$x\in
\Lambda_n$}
and every fixed $k\ge1$,
%
%e3.27 #&#
\begin{equation}
\label{eqi3} (1-\epsilon)c_{\infty}\leq\frac{\pb\{ A\cap\{
|S_k|>y, |S_i|\leq y, i\neq k\}\}}{n_1\pb\{ Y>x\} }\leq (1+
\epsilon)c_{\infty}.
\end{equation}
Write $z=x/(\log n)^{\xi}$ and introduce the probabilities, for $k\ge1$,
%
%e3.28 #&#
\begin{eqnarray}\label{eq09}
J_k&=&\pb \biggl\{A\cap\biggl\{ \biggl|\sum_{j\neq k}
(S_j- c_j)\biggr|>z,  |S_k|>y,
|S_i|\leq y, i\neq k\biggr\} \biggr\},
\nonumber\\
V_k&=& \pb \biggl\{A\cap\biggl\{ \biggl|\sum
_{j\neq k} (S_j- c_j)\biggr|\le z,
|S_k|>y, |S_i|\leq y, i\neq k\biggr\} \biggr\}.
\end{eqnarray}
Write $S=\sum(S_j-c_j)$, where summation is taken over the set
$\{j\dvtx  1\le j\le p, j\ne k,k\pm1,k\pm2\}$.
For $n$ sufficiently large, $J_k$ is
dominated by
\begin{eqnarray*}
&&\pb\bigl\{ |S|>z-8y, |S_k|>y, |S_i|\leq y, i\neq k\bigr\}\\
&&\qquad\leq
\pb\bigl\{ |S_k|>y\bigr\} \pb\bigl\{ |S|>0.5 z, |S_i|\leq y, i\neq k
\bigr\}.
\end{eqnarray*}
Applying the Prokhorov inequality (\ref{prokh}) in the same way as in
step 1b, we see that
\[
\pb\bigl\{ |S|>0.5 z, |S_i|\leq y, i\neq k\bigr\}\leq c n z^{-\a}
\leq c (\log n)^{-(M-\xi)\a}
\]
and by Markov's inequality,
\[
\pb\bigl\{ |S_1|>y\bigr\}\leq c\frac{n_1(y)}{y^{\a}} \leq c \frac{n_1}{y^{\a}}.
\]
Therefore
\[
\sup_{x\in\Lambda_n } \bigl(x^\a/n_1\bigr)
J_k\leq c (\log n)^{3\a\xi
-M\a}\to0.
\]
Thus it remains to bound the probabilities $V_k$. We start with
sandwich bounds
for~$V_k$,
%
%e3.29 #&#
%e3.30 #&#
\begin{eqnarray}
\label{eqvk}
&&\pb\bigl\{ S_k-c_k>x+z, |S|\leq z-8y,
|S_i|\leq y, i\neq k\bigr\}\nonumber\\[-8pt]\\[-8pt]
&&\qquad\le V_k\nonumber
\\
\label{eqvka}
&&\qquad\leq \pb\bigl\{ S_k-c_k>x-z, |S|\leq z+8y,
|S_i|\leq y, i\neq k\bigr\}.
\end{eqnarray}
By (\ref{es}), for every small $\epsilon>0$, $n$ sufficiently large
and uniformly for $x\in\Lambda_n$, we have
%
%e3.31 #&#
\begin{equation}
\label{eqineq}\qquad (1-\epsilon)c_{\infty} \leq\frac{\pb\{
S_k-c_k>x+z\}}{n_1 \pb\{ Y>x\}} \leq
\frac{\pb\{
S_k-c_k>x-z\}}{n_1 \pb\{ Y>x\}}\leq(1+\epsilon)c_{\infty},
\end{equation}
where we also used that $\lim_{n\to\infty}(x+z)/x=1$. Then the
following upper bound is immediate from (\ref{eqvka}):
\[
\frac{V_k}{n_1 \pb\{Y>y\}}\le\frac{\pb\{
S_k-c_k>x-z\}}{n_1 \pb\{ Y>x\}}\leq(1+\epsilon)c_{\infty}.
\]
From (\ref{eqvk}) we have
\[
\frac{V_k}{n_1 \pb
\{Y>y\}}\ge\frac{\pb\{ S_k-c_k>x+z\}}{n_1 \pb\{ Y>x\}}- L_k.
\]
In view of the lower bound in (\ref{eqineq}), the first
term on the right-hand side yields the desired lower bound if we
can show that {$L_k$} is negligible. Indeed, we have
\begin{eqnarray*}
n_1\pb\{ Y>x\} L_k&=&\pb\biggl\{
\{S_k-c_k>x+z\}\cap \biggl[\bigl\{|S|> z-8y\bigr\}\cup \bigcup
_{i\ne k}\bigl\{ |S_i|> y\bigr\} \biggr]\biggr\}
\\
&\le&\pb\bigl\{S_k-c_k>x+z,|S|> z-8y\bigr\}\\
&&{} + \sum
_{i\ne k}\pb\bigl\{S_k-c_k>x+z,
|S_i|> y\bigr\}
\\
&\le&\pb\{S_k-c_k>x+z\} \bigl[\pb\bigl\{|S|> z-8y\bigr\} +p \pb
\bigl\{|S_1|> y\bigr\} \bigr]
\\
&&{}+ c \sum_{|i-k|\le2,i\ne k} \pb\bigl\{S_k-c_k>x+z,
|S_i|> y\bigr\}.
\end{eqnarray*}
Similar bounds as in the proofs above yield that the
right-hand side is of the order $o(n_1/x^\a)$, hence $L_k=o(1)$. We
omit details. Thus we obtain (\ref{eqi3}).

\textit{Step} 1d: \textit{Final bounds.}
Now we are ready for the final steps in the proof
of (\ref{ss}) and (\ref{sss}).
Suppose {first} $c_{\infty}^+>0$ and
$\log x\le s_n$. In view of the decomposition (\ref{eqdecomp})
and steps 1a and 1b
we have as $\nto$ and
uniformly for $x\in\Lambda_n$,
\begin{eqnarray*}
&&
\frac{\pb\{
\sum_{j=1}^p (S_j-c_j)>x \}}{n \pb\{ Y>x\} }\\
&&\qquad\sim \frac{I_3}{n \pb
\{ Y>x\} }
\\
&&\qquad\sim\frac{n_1}{n}\frac{\sum_{k=1}^p\pb\{
\sum_{j=1}^p (S_j-\eb S_j)>x, |S_k|>y, |S_j|\leq y, j\neq
k\}}{n_1\pb\{ Y>x\} }.
\end{eqnarray*}
In view of step 1c, {in particular (\ref{eqi3}),}
the last expression
is dominated from above by
$(p n_1/n) (1+\epsilon)c_{\infty}\leq(1+\epsilon)c_{\infty}$ and from below
by
\[
\frac{n_1p}{n} (1-\epsilon)c_{\infty}\geq\frac{n-n_2-n_1} {
n}(1-
\epsilon)c_{\infty}\ge(1-\epsilon)c_{\infty} \biggl(1-
\frac{3s_n}{n\rho} \biggr).
\]
Letting first $\nto$ and then
$\epsilon\to0$, (\ref{sss}) follows and (\ref{ss}) is also satisfied
provided the additional condition
$\lim_{n\to\infty}s_n/n=0$ holds.

If $c_{\infty}^+=0$, then {$I_3=o(n \pb\{ |Y|>x\})$.} Let now $x\in
\Lambda_n$ and {recall the definition of $V_k$ from (\ref{eq09}).}
Then for every small $\d$ and sufficiently large $x$,
\begin{eqnarray*}
\frac{\pb\{ \sum_{j=1}^p (S_j-c_j)>x \}}{n \pb\{
|Y|>x\} }&\sim& \frac{I_3}{n \pb
\{ |Y|>x\} }
\\
&\le& \frac{n_1}{n}\frac{\sum_{k=1}^pV_k}{n_1\pb\{ |Y|>x\}}
\\
&\le& \sup_{x\in\Lambda_n}\frac{\pb\{ S_1>x(1-\d)-|c_1|\}} {
n_1\pb\{ |Y|>x\}},
\end{eqnarray*}
and (\ref{sss}) follows from the second
part of Lemma~\ref{estim}.

\textit{Step} 2: \textit{Proof of} (\ref{xx}) \textit{and} (\ref{zz}).
We restrict ourselves to
(\ref{xx}) since the proof of (\ref{zz}) is analogous.
Write $B=\{ | {\sum_{j=1}^{p_1}} (X_j-e_j)|>x\}$. Then
\[
\pb\{B\} \leq\pb \Biggl\{B \cap\bigcup_{k=1}^{p_1}
\bigl\{ |X_k|>y\bigr\} \Biggr\} + \pb \bigl\{ B\cap\bigl\{ |X_j|\leq y
\mbox{ for all $j\le p_1$}\bigr\} \bigr\}=P_1+P_2.
\]
By Lemma~\ref{x},
\[
P_1\le p_1 \pb\bigl\{ |X_1|>y\bigr\} \leq
C_1p_1y^{-\a}e^{-C_2(\log y)^{C_3}}=o
\bigl(nx^{-\a
}\bigr).
\]
For the estimation of $P_2$ consider
the random variables $X^y_j=X_j\one_{\{ |X_j|\leq y\} }$ and proceed
as in step 1b.\vspace*{9pt}

\textit{The case $\a>2$}.

The proof is analogous to
$\alpha\in(1,2]$. We indicate differences in
the main steps.

\textit{Step} 1b.
First we bound the large deviation probabilities of
the truncated sums (it is an analog of step 1b of Proposition~\ref{pomoc}).
We assume without loss of generality that $c_n\leq\log n$.
Our aim now is to prove that for $y=xc_n^{-0.5}$,
%
%e3.32 #&#
\begin{equation}
\label{eqg5} \quad\lim_{n\to\infty}\sup_{x\in\Lambda_n}\frac{x^\a}{n}\pb
\Biggl\{ \Biggl| \sum_{j=1}^p(S_j-
\eb S_j) \Biggr|>x, |S_j|\leq y \mbox{ for all $j\le p$}
\Biggr\} =0.
\end{equation}
We proceed as in the proof of (\ref{eqi1}) with
the same notation $S^y_j=S_j\one_{\{ |S_j|\leq y\}}$. As in the proof
mentioned,
$p |\eb S_j\one_{\{ |S_j|> y\} }|=o(x)$ and so we estimate the
probability of interest by
%
%e3.33 #&#
\begin{equation}
\label{eq300} I:= 3 \pb \biggl\{ \biggl|\sum_{1\le j\le p,j\in\{1,4,7,\ldots\}}
\bigl(S^y_j-\eb S^y_1\bigr) \biggr|
>x/6 \biggr\}.
\end{equation}
The summands in the latter sum are independent and therefore one can
apply the two-sided version of the Fuk--Nagaev inequality (\ref{fuknagaev})
to the random variables in
(\ref{eq300}): with $a_n=\b x/y=c_n^{0.5}\b$ and $p
\var(S^y_1)\le cpn_1^2\leq cnn_1$,
%
%e3.34 #&#
\begin{equation}
\label{eq88} I\le \biggl(c\frac{p n_1^{(1.5 +\sigma)\a}}{xy^{\a-1 }} \biggr)^{a_n} + \exp
\biggl\{-\frac{(1-\b)^2 x^2}{2e^{\a} cnn_1} \biggr\}.
\end{equation}
We suppose first that $x\leq n^{0.75}$. Then
the first quantity in (\ref{eq88}) multiplied by $x^{\a}/n$ is
dominated by
\begin{eqnarray*}
&&\bigl(c(\log n)^{(1.5 +\sigma)\a}c_n^{0.5(\a-1)}
\bigr)^{a_n} \bigl(n/x^\a\bigr)^{a_n-1} \\
&&\qquad\leq
c_n^{-0.5a_n(1+\a)+\a}\frac{(c(\log n)^{(1.5 +\sigma)\a
})^{a_n}}{(n^{0.5\a-1}(\log n)^{\a})^{a_n-1}}\to0.
\end{eqnarray*}
The second quantity
in (\ref{eq88}) multiplied by $x^{\a}/n$ is dominated by
\[
\frac{x^{\a}}{n} \exp \biggl\{-\frac{(1-\b)^2 c^2_n(\log n)^2} {
2e^{\a}cn_1} \biggr\}\leq
n^{\a\g-1}n^{-c c ^2_n}\to0.
\]
If
$x>n^{0.75}$ then $xn^{-0.5}>x^{\d}$ for an appropriately small
$\d$. Then the first quantity in (\ref{eq88}) is dominated by
\begin{eqnarray*}
&&
\bigl(c(\log x)^{(1.5 +\sigma)\a} \bigr)^{a_n}c_n^{0.5a_n(\a-1)}
\bigl(n/x^\a\bigr)^{a_n-1} \\
&&\qquad\leq c_n^{-0.5a_n(1+\a)+\a}
\frac{(c(\log x)^{(1.5 +\sigma)\a
})^{a_n}}{(n^{0.5\a-1}x^{\a\d})^{a_n-1}}\to0.
\end{eqnarray*}
The second
quantity is dominated by
\[
\frac{x^{\a}}{n} \exp \biggl\{-\frac{(1-\b)^2 c^2_nx^{2\d}(\log
n)^2} {2e^{\a}cn_1} \biggr\}\leq
x^{\a}e^{-cx^{\d}}\to0,
\]
which finishes the proof of (\ref{eqg5}).

\textit{Step} 1c. We prove
that, for any $\eps\in(0,1)$, sufficiently large $n$ and
fixed $k\ge1$, the following relation holds uniformly for $x\in
\Lambda_n$,
%
%e3.35 #&#
\begin{eqnarray}
\label{eqg7} (1-\eps)c_{\infty}&\leq&\frac{\pb\{ \sum_{j=1}^p (S_j-\eb S_j)>x,
|S_k|>y, |S_i|\leq y, i\neq k\} }{n_1 \pb\{ Y>x\} }\nonumber\\[-8pt]\\[-8pt]
&\leq& (1+
\eps)c_{\infty}.\nonumber
\end{eqnarray}
Let $z\in\Lambda_n$ be such that $x/z\to\infty$. As
for $\a\in(1,2]$, one
proves
\[
\frac{x^{\a}}{n_1}\pb \Biggl\{ \sum_{j=1}^p
(S_j-\eb S_j)>x, \biggl|\sum_{j\neq k}
(S_j-\eb S_j)\biggr|>z, |S_k|>y,
|S_j|\leq y, j\neq k \Biggr\} \to 0.
\]
Apply the Fuk--Nagaev inequality (\ref{fuknagaev})
to bound
\[
\pb \biggl\{ \biggl|\sum_{j\neq k} (S_j-\eb
S_j)\biggr|>\frac{z}{2}, |S_j|\leq y, j\neq k \biggr
\}.
\]
In the remainder of the proof one can follow the arguments of the proof
in step 1c for $\a\in(1,2]$.\vspace*{9pt}

\textit{The case $\a\leq1$}.
The proof is analogous to the case
$1< \a\leq2$; instead of Prokhorov's inequality (\ref{prokh})
we apply S. V. Nagaev's inequality (\ref{nagaev}).
We omit further details.
\end{pf}

%s3.5 #&#
\subsection{\texorpdfstring{Final steps in the proof of Theorem \protect\ref{mthm}}
{Final steps in the proof of Theorem 2.1}}
\label{subsecfinala}

We have for small $\vep>0 $,
%
%e3.36 #&#
\begin{eqnarray}
\label{eq77}
&&\pb\Biggl\{ \sum_{i=1}^n(
\wt S_i -\eb\wt S_i) >x(1+2\eps)\Biggr\}- \pb\Biggl\{
\biggl|\sum_{i=1}^n(\wt X_i -\eb
\wt X_i)\biggr| >x\eps\Biggr\} \nonumber\\
&&\quad{}-\pb\Biggl\{\Biggl| \sum
_{i=1}^n(\wt Z_i -\eb\wt
Z_i)\Biggr| >x\eps\Biggr\}
\nonumber
\\
&&\qquad\le \pb\{ \wt\S_n - \wt d_n>x\}
\\
&&\qquad \le \pb\Biggl\{ \sum_{i=1}^n(\wt
S_i -\eb\wt S_i) >x(1-2\eps)\Biggr\}+ \pb\Biggl\{
\sum_{i=1}^n(\wt X_i -\eb\wt
X_i) >x\eps\Biggr\} \nonumber\\
&&\qquad\quad{}+\pb\Biggl\{ \sum_{i=1}^n(
\wt Z_i -\eb\wt Z_i) >x\eps\Biggr\}.
\nonumber
\end{eqnarray}
Divide the last two probabilities in the first and last lines
by $n\pb\{ |Y|>x\}$. Then these ratios converge to
zero for $x\in\Lambda_n$, in view of (\ref{xx}), (\ref{zz}) and
Lemmas~\ref{x} and~\ref{z}. Now
\begin{eqnarray*}
&&
\pb\Biggl\{ \sum_{i=pn_1+1}^n(\wt
S_i -\eb\wt S_i) >x(1-2\eps)\Biggr\}\\
&&\qquad=\pb \Biggl\{
\sum_{i=pn_1+1}^{n-n_1}(\wt S_i -\eb
\wt S_i)>x(1-2\eps)\Biggr\}
\\
&&\qquad\leq\pb\Biggl\{ \sum_{i=pn_1+1}^{n-n_1}\bigl(
\underline{\wt S}_i>x(1-2\eps)\bigr) -\sum_{i=pn_1+1}^{n-n_1}|
\eb\wt S_i|\Biggr\},
\end{eqnarray*}
where $\underline{\wt S}_i=\Pi_{i,i+n_1-1}|B_{i+n_1}|+\cdots+\Pi_{i,i+n_2-1}|B_{i+n_2}|$.

Notice that if $i\leq n-n_2$ then (for $\a> 1$)
\[
\eb\wt S_i=\eb\wt S_1
\]
and for $n-n_2<i\leq n-n_1$
\[
\eb\wt S_i=(\eb A)^{n_1} \bigl(1+\eb A +\cdots+(\eb
A)^{n-i-n_1} \bigr)\eb B.
\]
Hence there is $C$ such that
\[
\sum_{i=pn_1+1}^{n-n_1}|\eb\wt S_i|
\leq2n_1C
\]
and so by
Proposition~\ref{block}
\begin{eqnarray*}
&&
\pb\Biggl\{ \sum_{i=pn_1+1}^{n-n_1}\bigl(\underline{
\wt S}_i>x(1-2\eps)\bigr) -\sum_{i=pn_1+1}^{n-n_1}|
\eb\wt S_i|\Biggr\}\\
&&\qquad\leq\pb\Biggl\{ \sum
_{i=1}^{n_1}\biggl(\underline{\wt S}_i>
\frac{x(1-2\eps) -2n_1C}{2}\biggr)\Biggr\}
\\
&&\qquad\quad{}+ \pb\Biggl\{ \sum_{i=n_1}^{2n_1}\biggl(\underline{
\wt S}_i>\frac{x(1-2\eps)
-2n_1C}{2}\biggr)\Biggr\}
\\
&&\qquad\leq C_1n_1x^{-\a}=o \bigl(n\pb\bigl\{ |Y|>x\bigr\}
\bigr),
\end{eqnarray*}
provided $\lim_{n\to\infty}s_n/n=0$.
Taking into account (\ref{ss}) and the sandwich
(\ref{eq77}), we conclude that
(\ref{res}) holds.

If the $x$-region is not bounded from above and $n> n_1(x)$ then
the above calculations together with Lemma~\ref{lemnew} show (\ref
{unres}). If $n\le
n_1(x)$, then
\[
\pb\{ \S_n- \wt d_n >x\} \leq C_1x^{-\a}e^{-C_2(\log
x)^{C_3}}
\]
and again (\ref{unres}) holds.

%s4 #&#
\section{Ruin probabilities}\label{ruin}
In this section we study
the \textit{ruin probability} related to the centered partial sum process
$T_n=\S_n-\eb\S_n$, $n\ge0$,
that is, for given $u>0$ and $\mu>0$ we consider the probability
\[
\psi(u)= \pb\Bigl\{\sup_{n \ge1} [T_n -\mu n ]>u\Bigr\}.
\]
We will work under the assumptions of Kesten's
Theorem~\ref{thmkesten}.
Therefore the random variables $Y_i$ are regularly varying with index
$\alpha>0$. Only
for
$\alpha>1$ the variable $Y$ has finite expectation and therefore we
will assume this condition throughout. Notice that the random walk
$(T_n-n \mu)$ has dependent steps and negative drift. Since
$(Y_n)$ is ergodic we have $n^{-1}(T_n-n \mu)\stas-\mu$ and in
particular $\sup_{n\ge1}(T_n-n \mu)<\infty$ a.s.

It is in general difficult to
calculate $\psi(u)$ for a given value $u$, and therefore most results
on ruin
study the asymptotic behavior of $\psi(u)$ when $u\to\infty$.
If the sequence $(Y_i)$ is i.i.d.\vadjust{\goodbreak} it is well known (see Embrechts and
Veraverbeke~\cite{embrechtsveraverbeke1982} for a special case of
subexponential step distribution and
Mikosch and Samorodnitsky~\cite{mikoschsamorodnitsky2000} for a
general regularly varying step distribution) that
%
%e4.1 #&#
\begin{equation}
\label{eqembrvera}
\psi_{\mathrm{ind}}(u) \sim\frac{u \pb\{Y>u\}}{\mu (\a-1)},\qquad u\to\infty.
\end{equation}
(We write $\psi_{\mathrm{ind}}$ to indicate that we are dealing with i.i.d.
steps.) If the step distribution has exponential moments the ruin probability
$\psi_{\mathrm{ind}}(u)$ decays to zero at an exponential rate; see, for
example,
the celebrated Cram\'er--Lundberg bound in Embrechts et al.
\cite{embrechtskluppelbergmikosch1997}, Chapter 2.

It is the main aim of this section to prove the following analog of
the classical ruin bound (\ref{eqembrvera}):
%
%th4.1 #&#
\begin{theorem}\label{mthm2} Assume that the conditions of
Theorem~\ref{thmkesten} are satisfied and additionally $B\ge0$ a.s.
and there exists $\eps>0$
such that $\eb A^{\a+\eps}$ and $\eb B^{\a+\eps}$ are finite, $\a
>1$ and $c_\infty^+>0$.
The following asymptotic result for the ruin probability holds for
fixed $\mu>0$, as $u\to\infty$:
%
%e4.2 #&#
\begin{eqnarray}
\label{eqruin} \P \Bigl\{ \sup_{n\ge1} (\S_n - \E
\S_n - n\mu) >u \Bigr\} &\sim& \frac{c_\infty}{\mu(\a-1)} u^{-\a+1}
\nonumber\\[-8pt]\\[-8pt]
&\sim&{ \frac{c_\infty}{c_\infty^+}\frac{u
\P\{Y>u\}}{\mu(\a-1)}.}\nonumber
\end{eqnarray}
\end{theorem}
%
%re4.2 #&#
\begin{remark}
We notice that the dependence in the sequence $(Y_t)$ manifests in the
constant $c_\infty/c_\infty^+$ in relation (\ref{eqruin}) which appears
in contrast to the i.i.d. case; see (\ref{eqembrvera}).
\end{remark}
To prove our result we proceed similarly as in the proof of
Theorem~\ref{mthm}. First notice that in view of (\ref{jm}),
\[
\P\Bigl\{ \sup_{n\ge1} \bigl(Y_0 \eta_n -
\E(Y_0 \eta_n) \bigr) >u \Bigr\} \le{\P\{
Y_0 \eta> u \} } = o\bigl(u^{1-\a}\bigr).
\]
Thus, it is sufficient to prove
\[
u^{\a-1} \P \Bigl\{ \sup_{n\ge1} (\wt\S_n - \E\wt
\S_n - n\mu) >u \Bigr\} \sim\frac{c_\infty}{\mu(\a-1)}
\]
for $\wt\S_n$ defined in (\ref{tsn}).
Next we change indices as indicated after (\ref{eq33}).
However, this time we cannot fix $n$ and therefore we will proceed
carefully; the details will be explained below.
Then we further decompose
$\wt\S_n$ into smaller pieces and study their asymptotic behavior.
%{\red(1+3\vep)^2} u^2 I_1
%&\le& \sum_{i=1}^{N-1} \eb(\wt R_i^2 \one_{D_i})\\
%&=& \sum_{i=1}^{N-1}\bigg( \sum_{j=1}^i \eb(\wt X_j^2\one_{D_i}) + 2
%&=& \sum_{i=1}^N\bigg( \sum_{j=1}^i \int_{A_i}\wt X_j^2 d\pb+ 2\sum_{1
%&\le& \sum_{i=1}^N\bigg( \sum_{j=1}^N \int_{A_i}\wt X_j^2 d\pb+ 2
%&\le& \E\big(\wt X_1 +\cdots+ \wt X_N\big)^2
%= \frac{N}{u^2} \eb\wt X_1^2.
%But I do not think we need Lemma 4.2.

\subsection*{Proof of Theorem \protect\ref{mthm2}}
The following lemma shows that the centered sums $(\wt\S_n-\eb\wt
S_n)_{n\ge
uM}$ for large $M$ do not contribute to the asymptotic behavior of the
ruin probability as $u\to\infty$.\vadjust{\goodbreak}
%
%le4.3 #&#
\begin{lemma}\label{lemoo} The following
relation holds:
\[
\lim_{M\to\infty} \limsup_{\uto} u^{\a-1} \pb \Bigl\{
\sup_{n>uM} (\wt\S_n -\eb\wt\S_n -n\mu)>u \Bigr
\}=0.
\]
\end{lemma}
\begin{pf}
Fix a large number $M$ and define the sequence
$N_l= uM 2^l$, $l\ge0$. Assume for the ease of
presentation that $(N_l)$ constitutes a sequence of even integers;
otherwise we can take $N_l=[uM]2^l$. Observe that
\[
\pb \Bigl\{ \sup_{n>uM} (\wt\S_n -\eb\wt\S_n -n
\mu)>u \Bigr\}\le \sum_{l=0}^\infty
p_l,
\]
where
$p_l=\pb \{
\max_{n\in[N_l,N_{l+1})} (\wt\S_n -\eb\wt\S_n
-n\mu)>u \}$.
For every fixed $l$, in the events above
we make the change of indices $i\to j=N_{l+1}-i$ and write, again
abusing notation,
\[
\wt Y_j= B_j +\Pi_{jj}B_{j+1}+
\cdots+\Pi_{j,N_{l+1}-2}B_{N_{l+1}-1}.
\]
With this notation, we have
\[
p_l=\pb \Biggl\{ \max_{n\in(0,N_l]} \sum
_{i=n}^{N_{l+1}-1} (\wt Y_i-\eb\wt
Y_i-\mu)>u \Biggr\}.
\]
Using the decomposition (\ref{eq2}) with the adjustment
$n_4=\min(j+n_3,\break N_{l+1}-1)$, we write $\wt Y_j=\wt U_j+\wt W_j$.
Then, by Lemma~\ref{red2}, for small $\delta>0$,
\begin{eqnarray*}
p_{l1}&=&\pb \Biggl\{ \max_{n\in
(0,N_l]} \sum
_{i=n}^{N_{l+1}-1}(\wt W_i - \eb\wt
W_i-\mu/4)>u/4 \Biggr\}
\\
&\le&\pb \Biggl\{ \sum_{i=N_l}^{N_{l+1}-1}(\wt
W_i-\eb\wt W_i)+ \sum_{i=1}^{N_{l}-1}
\wt W_i >u/4+N_l \mu/4 \Biggr\}
\\
&\le&\pb \Biggl\{ \sum_{i=1}^{N_{l+1}-1}(\wt
W_i-\eb\wt W_i) >u/4+N_l (\mu/4-E\wt
W_1) \Biggr\}
\\
&\le& c N_{l+1} N_l^{-\alpha-\delta}\le c
\bigl(uM2^l\bigr)^{1-\alpha-\delta}.
\end{eqnarray*}
We conclude that for every $M>0$,
\[
\sum_{l=0}^\infty p_{l1}= o
\bigl(u^{1-\alpha}\bigr) \qquad\mbox{as } u\to \infty.
\]
As in (\ref{eq4a}) we further decompose $\wt U_i=\wt X_i+\wt S_i+\wt
Z_i$, making the definitions precise in what follows.
Let $p$ be the smallest integer such that $pn_1\ge N_{l+1}-1$ for
$n_1=n_1(u)$.\vadjust{\goodbreak}
For $i=1,\ldots,p-1$ define $X_i$ as in (\ref{eq4}), and
$X_p=\sum_{i=(p-1)n_1+1}^{N_{l+1}-1}\wt X_i$. Now consider
\begin{eqnarray*}
p_{l2}&=&\pb \Biggl\{ \max_{n\in(0,N_l]} \sum
_{i=n}^{N_{l+1}-1}(\wt X_i - \eb\wt
X_i-\mu/4)>u/4 \Biggr\}
\\[-2pt]
&\le&\pb \Biggl\{ \sum_{i=N_l}^{N_{l+1}-1}(\wt
X_i - \eb\wt X_i)\\[-2pt]
&&\hspace*{12pt}{}+ \max_{n\in(0,N_l]} \sum
_{i=n}^{N_{l}-1}(\wt X_i - \eb\wt
X_i)>u/4+N_l \mu/4 \Biggr\}
\\[-2pt]
&\le&\pb \Biggl\{ \max_{k\le p/2}\sum_{i=k}^{p}(X_i
- \eb X_i) >u/8+N_l \mu/8 \Biggr\}
\\[-2pt]
&&{}+\pb \Biggl\{ \max_{k\le p/2}\sum_{i=k}^{p}(X_i
- \eb X_i) \le u/8+N_l \mu/8,\\[-2pt]
&&\hspace*{28pt} \max_{k\le p/2}
\max_{1\le j<n_1} \sum_{i=kn_1-j}^{kn_1}(\wt
X_i - \eb\wt X_i)>u/8+N_l \mu/8 \Biggr\}
\\[-2pt]
&=&p_{l21}+p_{l22}.
\end{eqnarray*}
The second quantity is estimated by
using Lemma~\ref{x} as follows for fixed $M>0$
\begin{eqnarray*}
p_{l22}&\le& c p \pb \Biggl\{ \sum_{i=1}^{n_1}
\wt X_i>u/8+N_l \mu/8 \Biggr\} \le C_1 p
N_l^{-\a}e^{-C_2 (\log N_l)^{C_3}}\\[-2pt]
&=&o\bigl(u^{1-\a}\bigr)
2^{-{
(\a-1)} l},
\end{eqnarray*}
where $C_i$, $i=1,2,3$, are some positive
constants. %, $\epsilon>0$ is chosen sufficiently small.
Therefore
for every fixed~$M$,
\[
\sum_{l=0}^\infty p_{l22}= o
\bigl(u^{1-\alpha}\bigr)\qquad \mbox{as } u\to\infty.
\]
Next we treat
$p_{l21}$. We observe that $X_i$ and $X_j$ are independent for
$|i-j|>1$. Splitting summation in $p_{l21}$ into the subsets of
even and odd integers, we obtain an estimate of the following type
\[
p_{l21}\le c_1\pb \biggl\{ \max_{k\le p/2}\sum
_{k\le2i\le
p}(X_{2i} - \eb X_{2i})
>c_2 (u+N_l) \biggr\},
\]
where the summands are now independent. By the law of large numbers,
for any $\epsilon\in(0,1)$, $k\le p/2$, large $l$,
\[
\pb \biggl\{ \sum_{2i\le k}(X_{2i} - \eb
X_{2i}) > -\epsilon c_2 (u+N_l) \biggr\}\ge1/2.\vadjust{\goodbreak}
\]
An application of the maximal inequality (\ref{levy}) in the \hyperref[app]{Appendix}
and an adaptation of Proposition~\ref{pomoc}
yield
\[
p_{l21}\le2 \pb \biggl\{ \sum_{2i\le p}(X_{2i}
- \eb X_{2i}) >(1-\epsilon)c_2 (u+N_l) \biggr
\} \le c N_l^{1-\alpha}.
\]
Using the latter bound and summarizing the above estimates, we finally
proved that
\[
\lim_{M\to\infty}\limsup_{u\to\infty} u^{\alpha-1} \sum
_{l=0}^\infty p_{l2}=0.
\]
Similar arguments show that the sums involving the $\wt S_i$'s and $\wt
Z_i$'s are negligible as well. This proves the lemma.
\end{pf}
In view of Lemma~\ref{lemoo} it suffices to study the following probabilities
for sufficiently large $M>0$:
\[
\P \Bigl\{ \max_{n\le Mu} (\wt\S_n - \E\wt\S_n - n
\mu) > u \Bigr\}.
\]
Write $N_0 = \lfloor Mu \rfloor$, change again indices
$i\to j = N_0 -i+1$ and write, abusing notation,
\[
\wt Y_j = B_j + \Pi_{jj}B_{j+1} +
\cdots+ \Pi_{j,N_0-1} B_{N_0}.
\]
Then we decompose $\wt Y_j$ as in the proof of Lemma~\ref{lemoo}.
Reasoning in the same way as above, one proves that the
probabilities related to the quantities $\wt W_i$, $\wt X_i$ and
$\wt Z_i$ are of lower order than $u^{1-\alpha}$ as $u\to\infty$
and, thus, it remains to study the probabilities
%
%e4.3 #&#
\begin{equation}
\label{eqglow} \P \Biggl\{ \max_{n\le
N_0} \sum
_{i=n}^{N_0} (\wt S_i - \E\wt
S_i - \mu) > u \Biggr\},
\end{equation}
where $\wt S_i$ were defined in (\ref{eq4a}).

Take $n_1 = \lfloor\log N_0/\rho\rfloor$,
$p=\lfloor N_0/n_1\rfloor$ and denote by $S_i$ the sums of $n_1$
consecutive $\wt S_i$'s as
defined in (\ref{eq4}).
Observe that for any $n$ such that
$n_1(k-1)+1\leq n \leq kn_1$, $k-1\leq p$ we have
\begin{eqnarray*}
\sum_{i=n}^{N_0} (\wt S_i -
\E\wt S_i - \mu)&\leq&2n_1 (\E\wt S_1 +
\mu)+\sum_{i=(k-1)n_1+1}^{(p+1)n_1}(\wt S_i
- \E\wt S_i - \mu)
\\
&\leq&2n_1 (\E\wt S_1 +\mu)+\sum
_{i=k-1}^{p}( S_i - \E S_i
- n_1\mu)
\end{eqnarray*}
and
\begin{eqnarray*}
\sum_{i=n}^{N_0} (\wt S_i -
\E\wt S_i - \mu)&\geq&-2n_1 (\E\wt S_1 +
\mu)+\sum_{i=kn_1}^{pn_1}(\wt S_i
- \E\wt S_i - \mu)
\\[-2pt]
&\geq&-2n_1 (\E\wt S_1 +\mu)+\sum
_{i=k}^{p-1}( S_i - \E
S_i - n_1\mu).
\end{eqnarray*}
Therefore and since $n_1$ is of order $\log u $,
instead of the probabilities (\ref{eqglow}) it suffices to study
\[
\psi_p(u)=\P \Biggl\{ \max_{n\le p} \sum
_{i=n}^p (S_i - \E S_i
- n_1\mu) > u \Biggr\}.
\]
%
%where $p=\lfloor N_0/n_1\rfloor$, $n_1 = \lfloor\log N_0/\rho
%defined in \eqref{eq4}.

Choose $q =[M/ \eps_1^{\a+1}]+1$ for some small
$\vep_1$ and large $M$. Then the random variables
\[
R_k = \sum_{i=(k-1)q}^{kq-3}
S_i,\qquad k=1,\ldots,r = \lfloor p/q\rfloor,
\]
are independent and we have
\begin{eqnarray*}
\psi_p(u)&\le& \P \biggl\{ \max_{n\le qr} \mathop{\sum
_{ n\le i\le qr}}_{i\not= kq-2,kq-1} (S_i - \E
S_i - n_1\mu) > u(1-3\eps_1) \biggr\}
\\[-2pt]
&&{}+ \P \Biggl\{ \max_{j\le r} \sum_{k=j}^r
(S_{kq-2} - \E S_{kq-2} - n_1\mu) >
\eps_1 u \Biggr\}
\\[-2pt]
&&{}+ \P \Biggl\{ \max_{j\le r} \sum_{k=j}^r
(S_{kq-1} - \E S_{kq-1} - n_1\mu) >
\eps_1 u \Biggr\}
\\[-2pt]
&&{}+ \P \Biggl\{ \max_{qr<n<p} \sum_{i=n}^p(S_i
- \E S_i -n_1\mu) > \eps_1 u \Biggr\}
=: \sum_{i=1}^4\psi_{p}^{(i)}(u).
\end{eqnarray*}
The quantities $\psi_p^{(i)}(u)$, $i=2,3$, can be estimated in the
same way; we focus on $\psi_p^{(2)}(u)$.
Applying Petrov's inequality (\ref{petrov}) and Proposition~\ref{pomoc},
we obtain for some constant $c$ not depending on $\vep_1$,
\begin{eqnarray*}
\psi_p^{(2)}(u) &\le& \P \Biggl\{ \max_{j\le r}
\sum_{k=j}^r (S_{kq-2} - \E
S_{kq-2}) > \eps_1 u \Biggr\}
\\[-2pt]
&\le& c \P \Biggl\{ \sum_{k=1}^r
(S_{kq-2} - \E S_{kq-2}) >\eps_1 u/2 \Biggr\}
\\[-2pt]
&\le& c r n_1 (\eps_1 u)^{-\a} \le c
\eps_1 u^{-\a+1}.
\end{eqnarray*}
Hence we obtain $\lim_{\vep_1\downarrow0} \limsup_{u\to\infty}
u^{\a-1} \psi_p^{(2)}(u)=0$. By (\ref{sss}), for some constant $c$,
\[
\psi_p^{(4)}(u)\leq c\frac{qn_1}{(\eps_1u)^{\a}}\leq c
\frac{Mu}{r(\eps_1u)^{\a}}.
\]
Since $r\geq q>
M/ \eps_1^{\a+1}$
for large $u$ we
conclude for such $u$ that
$r^{-1}\leq M^{-1}\eps_1^{\a+1}$ and therefore
$
\psi_p^{(4)}(u)\leq c\eps_1u^{1-\a}$ and
$\lim_{\vep_1\downarrow0} \limsup_{u\to\infty}
u^{\a-1} \psi_p^{(4)}(u)=0$.

Since $A$ and $B$ are nonnegative we have for large $u$ with
$\mu_0 = \mu(q-2)$,
\begin{eqnarray*}
\psi_p^{(1)}(u) &\le& \P \Biggl\{ \max_{j\le r}
\sum_{i=1}^j (R_i - \E
R_i - \mu_0 n_1) > u (1-3
\eps_1) - qn_1(\eb S_1+\mu) \Biggr\}
\\
&\le& \P \Biggl\{ \max_{j\le r} \sum_{i=1}^j
(R_i - \E R_i - \mu_0
n_1) > u (1-4\eps_1) \Biggr\}.
\end{eqnarray*}
Combining the bounds above we proved
for large $u$, small $\vep_1$ and some constant $c>0$ independent
of $\vep_1$ that
\[
\psi_p(u) \le \P \Biggl\{ \max_{{ j}\le r} \sum
_{i=1}^{{ j}} (R_i - \E
R_i - \mu_0 n_1) > u (1-4
\eps_1) \Biggr\} + c \eps_1 u^{-\a+1}.
\]
Similar arguments as above show that
\[
\psi_p(u) \ge \P \Biggl\{ \max_{{ j}\le r} \sum
_{i=1}^{{ j}} (R_i - \E
R_i - \mu_0 n_1) > u (1+4
\eps_1) \Biggr\} - c \eps_1 u^{-\a+1}.
\]
Thus we reduced the problem to study an expression consisting of
independent random variables $R_i$ and
the proof of the theorem is finished if we can show the
following result. Write
\[
\Omega_r= \Biggl\{ \max_{j\le r } \sum
_{i=1}^j (R_i-\E
R_i-n_1\mu_0)>u \Biggr\}.
\]

%le4.4 #&#
\begin{lemma}
The following relation holds
\[
\lim_{M\to\infty}\limsup_{\uto} \biggl|{u^{\a-1}}\P\{
\Omega_r \} - \frac{c_\infty c_\infty^+}{\mu(\a-1)} \biggr|=0.
\]
\end{lemma}
%
%Let $u_0 = u (1-4\eps_1)$, $\mu_0 = \mu(q-2)$,
%$C_0 = (q-2) (c_\infty^+ c_\infty)/(c_\infty^+ + c_\infty^-)$. Then
%for every small $\eps_1$ one can find sufficiently large $M$
%such that for large $u_0$
%ele
%
\begin{pf}
Fix some $\eps_0>0$ and choose some large $M$. Define $C_0= (q-2)
c_\infty
c_\infty^+$.
Reasoning as in the proof of (\ref{ss}), we obtain
for any integers $0\le j < k \le r$ and large $u$
%
%e4.4 #&#
\begin{equation}
\label{eqRjk} 1-\eps_0\le u^\a \frac{\pb \{ \sum_{i=j+1}^k  ( R_i - \E R_i  )
>u  \}}{
n_1(k-j)C_0}
\le1+\vep_0.
\end{equation}
%
%Moreover, since we consider only positive $B_i$'s
Choose $\eps, \d>0$ small to be determined later in dependence on
$\vep_0$. Eventually, both $\eps, \d>0$ will\vspace*{1pt} become arbitrarily\vadjust{\goodbreak} small
when $\vep_0$ converges to zero. Define the
sequence $k_0 = 0$, $k_l = [\d n_1^{-1}(1+\eps)^{l-1}u]$, $l\ge1$.
Without loss of generality we will assume $k_{l_0}= Mu n_1^{-1}$
for some integer number $l_0$.
%We consider the block sums
%$\sum_{i=k_l+1}^{k_{l+1}}(R_i-\E R_i - \mu_0 n_1)$ containing
%$k_{l+1} - k_l = \d\eps n_1^{-1} (1+\eps)^{l-1}u$ summands.
For $\eta> \eps_0 (2l_0)^{-1}$ consider the independent events
\[
D_l = \Biggl\{ \max_{k_l < j \le k_{l+1}} \sum
_{i=k_l+1}^j (R_i - \E
R_i) > 2\eta u \Biggr\},\qquad l=0,\ldots, l_0-1.
\]
Define the disjoint sets
\[
W_l =\Omega_r \cap D_l \cap\bigcap
_{m\not= l} D_m^c,\qquad l=0,\ldots,l_0-1.
\]
We will show that
%
%e4.5 #&#
\begin{equation}
\label{eqqq} \Biggl| \P\{ \Omega_r\} - \sum_{l=0}^{l_0-1}
\P\{W_l\} \Biggr| \le o\bigl(u^{1-\a}\bigr),\qquad \uto.
\end{equation}
First we observe that
$
\Omega_r\subset\bigcup_{l=0}^{l_0-1}D_l.
$
Indeed, on $\bigcap_{l=0}^{l_0-1}D_l^c$ we have
\[
\max_{j\le r } \sum_{i=1}^j
(R_i-\E R_i-n_1\mu_0)\le
l_0 2\eta u\le\vep_0 u,
\]
and therefore $\Omega_r$ cannot hold for small $\vep_0$.
Next we prove that
%
%e4.6 #&#
\begin{equation}
\label{eqdl} \P\biggl\{ \bigcup_{m\not=l}
(D_m\cap D_l) \biggr\} = o\bigl(u^{1-\a}\bigr),\qquad \uto.
\end{equation}
Then (\ref{eqqq}) will follow. First apply Petrov's inequality
(\ref{petrov}) to $\P\{D_l\}$ with $q_0$ arbitrarily close to one and
power $p_0\in(1,\a)$. Notice that $\E|R_i|^{p_0}$ is of the
order~$qn_1$, hence $m_{p_0}$ is of the order
%{ Here one needs an upper bound for $q$ as well}
$\d\eps qu\le c \d\eps
M\eps_1^{-\a-1}u$. Next apply
(\ref{eqRjk}). Then one obtains for sufficiently large $u$, and
small $\vep,\delta$, and some constant $c'$ depending on
$\vep,\delta,\vep_0,\vep_1$,
\begin{eqnarray*}
\P\{D_l\} &\le&q_0^{-1}\P \Biggl\{ \sum
_{i=k_l+1}^{k_{l+1}} (R_i -\E
R_i) > \eta u \Biggr\}
\\
&\le&q_0^{-1}n_1(k_{l+1}-k_l)
(1+\eps_0)C_0 (\eta u)^{-\a}\le
c'u^{1-\a}.
\end{eqnarray*}
Hence
$
\P\{\bigcup_{m\neq l}(D_l\cap D_m)\}=O(u^{2(1-\alpha)})$
as
desired for (\ref{eqdl}) if all
the parameters $\vep,\delta,\vep_0,\vep_1$ are fixed.
%%since we may choose sufficiently small $\eps$.} Notice that $$
%$$
%Indeed on the set $\bigcap_{l\le l_0} D_l^c$ we have
%$$
%$$
%Thus, we are led to study
%probabilities of the sets
%$$
%W_l = \bigg\{ \max_{j\le r}\sum_{i=1}^j (R_i - \E R_i -\mu_0 n_1) > u
%$$ since $\P\{ \bigcup_{m\not=l} (D_m\cap D_l) \} = o(u^{2(1-\a)})$ and
%$$
%- \sum_{l=0}^{l_0-1}\P\{W_l\}
%$$

Thus we showed (\ref{eqqq}) and it remains to find suitable bounds
for the probabilities $\pb\{W_l\}$.
On the set $W_l$ we have
\begin{eqnarray*}
\max_{j\le k_l} \sum_{i=1}^j
(R_i - \E R_i-\mu_0 n_1)&
\le& \max_{j\le k_l} \sum_{i=1}^j
(R_i - \E R_i) \le2\eta lu \le\eps_0 u,
\\
\max_{k_{l+1} < j \le r}\sum_{i=k_l+1}^j(R_i
- \E R_i)&\le& 2\eta l_0u \le\eps_0 u.
\end{eqnarray*}
Therefore for small $\vep_0$ and large $u$ on the event $W_l$,
\begin{eqnarray*}
&&
\max_{j \le r} \sum_{i=1}^j(R_i
-\E R_i - \mu_0 n_1) \\[-2pt]
&&\qquad=
\max_{k_l< j\le r} \sum_{i=1}^j
(R_i -\E R_i - \mu_0n_1)
\\[-2pt]
&&\qquad\le \sum_{i=1}^{k_l}(R_i -
\E R_i - \mu_0 n_1) +
\max_{k_l < j \le
k_{l+1}}\sum_{i=k_l+1}^j
(R_i - \E R_i) \\[-2pt]
&&\qquad\quad{}+ \max_{k_{l+1} < j \le r} \sum
_{i=k_{l+1}+1}^j (R_i -\E R_i)
\\[-2pt]
&&\qquad\le { 2}\eps_0 u - k_l \mu_0
n_1 + \max_{k_l < j \le k_{l+1}}\sum_{i=k_l+1}^j(R_i
- \E R_i).
\end{eqnarray*}
Petrov's inequality (\ref{petrov}) and (\ref{eqRjk}) imply
for
$l\ge1$ and large $u$ that
\begin{eqnarray*}
\P\{W_l\}&\le& \P \Biggl\{ \max_{k_l < j \le k_{l+1}} \sum
_{i=k_l+1}^j (R_i-\E R_i)
\ge(1-2\eps_0)u + \mu_0 n_1
k_l \Biggr\}
\\[-2pt]
&\le& q_0^{-1} \P \Biggl\{ \sum
_{i=k_l+1}^{k_{l+1}} (R_i-\E R_i)
\ge (1-3\eps_0)u + \mu_0 n_1
k_l \Biggr\}
\\[-2pt]
&\le& q_0^{-1} \frac{ (k_{l+1}-k_l)n_1(1+\eps_0) C_0}{((1-3\eps_0)
+ \mu_0\d(1+\eps)^{l-1})^\a u^\a}
\\[-2pt]
&=& q_0^{-1}\frac{ \d\eps(1+\eps)^{l-1} (1+\eps_0)C_0 }{((1-{
3}\eps_0) + \mu_0\d(1+\eps)^{l-1})^\a} u^{1-\a}.
\end{eqnarray*}
For $l=0$ we use exactly the same arguments, but in this case
$(k_1-k_0)n_1 = \d u$ and $k_0 = 0$. Thus we arrive at the upper bound
%
%e4.7 #&#
\begin{eqnarray}
\label{eqsumWa}\quad \sum_{l=0}^{l_0-1} \P
\{W_l\} &\le& q_0^{-1}(1+\eps_0)\nonumber\\[-2pt]
&&{}\times
C_0 \Biggl( \frac{\d}{(1-{ 3}\eps_0)^\a} + \sum_{i=1}^{l_0-1}
\frac{
\d\eps(1+\eps)^{l-1}} {((1-{ 3}\eps_0) + \mu_0\d(1+\eps
)^{l-1})^\a} \Biggr) u^{1-\a}
\\[-2pt]
&=& q_0^{-1}(1+\vep_0) A(\eps, \d,
\eps_0,l_0)u^{1-\a}.
\nonumber
\end{eqnarray}
To estimate $\P\{W_l\}$ from below first notice that on $W_l$, for
large $u$,
\begin{eqnarray*}
\max_{j\le r}\sum_{i=1}^j
(R_i - \E R_i -\mu_0 n_1)
&\ge& \sum_{i=1}^{k_{l+1}} (R_i -
\E R_i -\mu_0 n_1)
\\
&\ge& \sum_{i= k_l+1}^{k_{l+1}} (R_i
- \E R_i) - k_{l+1}\mu_0 n_1
- k_l\E R_1
\\
&\ge& \sum_{i= k_l+1}^{k_{l+1}} (R_i
- \E R_i) - k_{l+1}\mu_0 n_1
- \eps_0 u.
\end{eqnarray*}
By (\ref{eqdl}) and (\ref{eqRjk}), we have for $l\ge 1$ and as
$u\to\infty$,
\begin{eqnarray*}
\P\{W_l\} &\ge& \P \Biggl\{ \Biggl\{ \sum
_{i=k_l+1}^{k_{l+1}} (R_i-\E
R_i)> (1+\eps_0)u + \mu_0
n_1 k_{l+1} \Biggr\} \cap D_l\cap\bigcap
_{m\not=l}D_m^c \Biggr\}
\\
&\ge& \P \Biggl\{ \sum_{i=k_l}^{k_{l+1}-1}
(R_i-\E R_i)> (1+\eps_0)u +
\mu_0 n_1 k_{l+1} \Biggr\} -\P \biggl\{
D_l\cap\bigcup_{m\not
=l}D_m
\biggr\}
\\
% &&- \P\bigg\{ \bigg\{ \sum_{i=k_l}^{k_{l+1}-1} (R_i-\E R_i)> (1+
&\ge& \frac{(k_{l+1}-k_l)n_1(1-\eps_0)C_0}{( (1+\eps_0)u+\mu_0
k_{l+1}n_1 )^\a} - { o
\bigl(u^{1-\a}\bigr)} \ge\frac{(1-{ 2}\eps_0) C_0 \d(1+\eps)^{l-1} \eps}{(
(1+\eps_0)+\mu_0\d(1+\eps)^l)^\a} u^{1-\a}.
\end{eqnarray*}
Hence
\begin{eqnarray*}
\sum_{l=0}^{l_0-1} \P\{W_l\} &
\ge&(1-{ 2}\eps_0)\\
&&{}\times C_0 \Biggl(\frac{\d}{(1+\eps_0+\mu_0 \d)^\a} +\sum
_{l=1}^{l_0-1} \frac{\d(1+\eps)^{l-1} \eps}{( (1+\eps_0)+\mu_0\d(1+\eps
)^l)^\a} \Biggr)
u^{1-\a}
\\
& =&(1-{ 2}\eps_0)C_0 B(\eps, \d,
\eps_0,l_0)u^{1-\a}.
\end{eqnarray*}
Thus we proved that
%
%e4.8 #&#
\begin{eqnarray}
\label{eq00b}
(1-2\vep_0)B(\eps, \d,
\eps_0,l_0) &\leq& \liminf_{u\to\infty}
C_0^{-1}u^{\a-1}\sum
_{l=0}^{l_0-1} \P\{W_l\}
\nonumber
\\
& \leq& \limsup_{u\to\infty}C_0^{-1} u^{\a
-1}
\sum_{l=0}^{l_0-1} \P\{W_l\}\\
&\leq&
q_0^{-1} (1+\vep_0) A(\eps, \d,
\eps_0,l_0).\nonumber
\end{eqnarray}
Finally, we will justify that the upper and lower bounds are
close for small $\vep,\delta,\vep_0$, large $M$ and $q_0$ close to 1.
For a real number $s$ which is small in absolute value
define the functions
\[
f_s(x) = (1+s+\mu_0x)^{-\a} \quad\mbox{and}\quad
F_s(x) = (1+s+\mu_0 x) f_s(x) \qquad\mbox{on
$[\d, M]$}.
\]
Let $x_l = \d
(1+\eps)^{l-1}$, $l=1,\ldots,l_0$. Since $x_{l+1}-x_l = \d\eps
(1+\eps)^{l-1}$ are uniformly bounded by $\eps M$ and $f_s$ is
Riemann integrable on $[0,\infty)$, choosing $\eps$ small, we have
\begin{eqnarray*}
A(\vep,\d,\vep_0,l_0) &=& \sum
_{l=1}^{l_0-1} f_{-3\eps_0}(x_l)
(x_{l+1}-x_l)
\\
&\le& \int_\d^M f_{- 3 \eps_0}(x) \,dx =
\frac{F_{-3\eps_0}(\d) -
F_{-3\eps_0}(M)+\vep_0}{\mu_0 (\a-1)}.
\end{eqnarray*}
Thus we obtain the bound
%
%e4.9 #&#
\begin{equation}
\label{eq01a} \lim_{q_0\uparrow1} \lim_{\vep_0\downarrow0}\lim_{M\to\infty}
\lim_{\delta\downarrow
0}q_0^{-1} (1+\vep_0) A(
\eps, \d, \eps_0,l_0)= \bigl(\mu_0 (
\a-1)\bigr)^{-1}.
\end{equation}
Proceeding in a similar way,
\[
B(\eps, \d, \eps_0,l_0) \ge\int_\d^M
f_{ \eps_0}(x)\,dx = \frac{F_{\eps_0}(\d) -
F_{\eps_0}(M)-\vep_0}{\mu_0 (\a-1)}.
\]
The right-hand side converges to $(\mu_0 (\a-1))^{-1}$ by letting
$\delta\downarrow0$, $M\to\infty$ and $\vep_0\downarrow0$.
The latter limit relation in combination with
(\ref{eq00b}) and (\ref{eq01a})
proves the lemma.
% \\ \le
%u^{\a-1} \P\Big\{ \bigcup_{l=0}^{l_0-1} W_l \Big\} \le
%(1+\eps_0)C_0\bigg( \frac{F_{-2\eps_0}(0)-F_{-2\eps_0}(M)}{\mu_0(
\end{pf}
%
%sA #&#
\begin{appendix}\label{app}
%sA #&#
\section{Inequalities for sums of independent random~variables}
In this section, we consider a sequence $(X_n)$ of independent
random variables and their partial
sums $R_n=X_1+\cdots+X_n$. We always write
$B_n=\var(R_n)$ and $m_{p }=\sum_{j=1}^n\eb|X_j|^{p } $ for $p>0$.
First, we collect some of the classical
tail estimates for $R_n$.
%
%leA.1 #&#
\begin{lemma} The following inequalities hold.

\textup{Prokhorov's inequality} (cf. Petrov~\cite{petrov1995},
page 77): Assume that the $X_n$'s are centered, $|X_n|\leq y$ for all
$n\ge1$ and some $y>0$. Then
%
%eA.1 #&#
\begin{equation}
\label{prokh} \pb\{ R_n\geq x\} \leq\exp \biggl\{ -
\frac{x}{2 y}\operatorname{arsinh} \biggl(\frac{xy}{2 B_n} \biggr) \biggr\},\qquad x>0.
\end{equation}

\textup{S. V.
Nagaev's inequality} (see~\cite{nagaev1979}): Assume $m_{p }<\infty$
for some $p>0$. Then
%
%eA.2 #&#
\begin{equation}
\label{nagaev} \pb\{R_n>x\} \leq\sum_{j=1}^n
\pb\{ X_j>y\} + \biggl(\frac{e m_{p }}{xy^{p -1}} \biggr)^{x/y},\qquad
x,y>0.
\end{equation}

\textup{Fuk--Nagaev inequality} (cf.
Petrov~\cite{petrov1995}, page 78): Assume that the $X_n$'s are
centered, $p >2$, $\b=p /(p +2)$ and $m_{p }<\infty$. Then
%
%eA.3 #&#
\begin{eqnarray}
\label{fuknagaev} \pb\{ R_n>x\} &\leq&\sum
_{j=1}^n\pb\{ X_j>y\}+ \biggl(
\frac{m_{p }}{\b x y^{p -1}} \biggr)^{\b x/y}\nonumber\\[-8pt]\\[-8pt]
&&{}+ \exp \biggl\{ -\frac{(1-\b)^2x^2}{2e^{p }B_n}
\biggr\},\qquad x,y>0.\nonumber
\end{eqnarray}

\textup{Petrov's inequality} (cf. Petrov~\cite{petrov1995},
page 81): Assume that the $X_n$'s are centered and $m_p<\infty$ for
some $p\in(1,2]$. Then for every $q_0<1$, with $L=1$ for $p=2$
and $L=2$ for $p\in(1,2)$,
%
%eA.4 #&#
\begin{equation}
\label{petrov} \pb\Bigl\{\max_{i\le
n} R_i>x\Bigr\} \leq
q_0^{-1} \pb\bigl\{R_n>x - \bigl[
\bigl(L/(1-q_0)\bigr)^{-1}m_p
\bigr]^{1/p}\bigr\},\qquad x\in\bbr.\hspace*{-30pt}
\end{equation}

\textup{L\'evy--Ottaviani--Skorokhod inequality} (cf.
Petrov~\cite{petrov1995}, Theorem 2.3 on pa\-ge~51): If
$\pb\{R_n-R_k\ge-c)\}\ge q$, $k=1,\ldots,n-1$, for some constants
$c\ge0$ and $q>0$, then
%
%eA.5 #&#
\begin{equation}
\label{levy} \pb\Bigl\{\max_{i\le n} R_i>x\Bigr\} \leq
q^{-1} \pb\{R_n>x -c\},\qquad x\in\bbr.
\end{equation}
\end{lemma}

%sA #&#
\section{\texorpdfstring{Proof of Lemma \lowercase{\protect\ref{estim}}}{Proof of Lemma 3.7}}\label{app2}

Assume first $c_\infty^+>0$. We have by independence of $Y$ and $\eta_k$,
for any
$k\ge1$, $x>0$ and $r>0$,
\[
\frac{\pb\{ \eta_kY>x\}}{k \pb\{ Y>x\} } = \biggl(\int_{(0,x/r]}+\int
_{[x/r,\infty)} \biggr) \frac{\pb\{
Y>x/z\}}{k \pb\{ Y>x\} } \,d\pb(\eta_k
\leq z) =I_1+I_2.
\]
For every $\eps\in(0,1)$ there is $r>0$ such that for $x\ge r$ and
$z\le x/r$,
\[
\frac{\pb\{ Y>x/z \}}{\pb\{ Y>x\} }\in z^{\a}[1-\eps, 1+\eps ]
\quad\mbox{and}\quad \pb\{ Y>x
\}x^{\a}\geq c_\infty^+-\eps.
\]
Hence for sufficiently large $x$,
\[
I_1\in k^{-1}\eb\eta_k^{\a}
\one_{\{\eta_k\leq x/r\} }[1-\eps,1+\eps]
\]
and
\[
I_2\leq c
k^{-1}x^{\alpha}\pb\{\eta_k>x/r\}\le c
k^{-1} \eb \eta_k^\alpha\one_{\{\eta_k>x/r\}}.
\]
We have
\[
I_1\in\bigl(k^{-1}\eb\eta_k^{\a}-k^{-1}
\eb\eta_k^{\a}\one_{\{\eta_k>
x/r\}}\bigr)[1-\eps,1+\eps]
\]
and by virtue of Bartkiewicz et al.
\cite{bartkiewiczjakubowskimikoschwintenberger2009},
$\lim_{k\to\infty}k^{-1}\eb\eta_k^{\a} = c_\infty$.
Therefore it is enough to prove that
%
%eA.1 #&#
\begin{equation}
\label{eqnull} \lim_{n\to\infty} \sup_{r_n\le k\le n_1, b_n\le x} k^{-1}\eb
\eta_k^{\a}\one_{\{\eta_k>
x\}}=0.
\end{equation}
By the H\"older and Markov inequalities
we have for $\epsilon>0 $,
%
%eA.2 #&#
\begin{equation}
\label{eq6} \eb\eta_k^{\a} \one_{\{\eta_k> x\}}\leq\bigl(
\eb\eta_k^{\alpha+\epsilon}\bigr)^
{\a/(\a+\epsilon) } \bigl(\pb\{
\eta_k>x\} \bigr)^{\epsilon/(\alpha+\epsilon) } \leq x^{-\epsilon} \eb
\eta_k^{\a+\epsilon}.
\end{equation}
Next we study the order of magnitude of $\eb\eta_k^{\a+\epsilon}$.
By definition of $\eta_k$,
\begin{eqnarray*}
\eb\eta_k^{\a+\epsilon}&=&\eb A^{\a+\epsilon}\eb(1+
\eta_{k-1})^{\a+\epsilon}
\\
&=&\eb A^{\a+\epsilon} \bigl(\eb(1+\eta_{k-1})^{\a+\epsilon}-\eb
\bigl(\eta_{k-1}^{\a+\epsilon
}\bigr) \bigr) + \eb A^{\a+\epsilon}\eb
\eta_{k-1}^{\a+\epsilon}.
\end{eqnarray*}
Thus we get the recursive relation
%
%eA.3 #&#
\begin{eqnarray}\label{eq7}
\eb\eta_k^{\a+\epsilon}&=&\sum_{i=1}^k
\bigl(\eb A^{\a
+\epsilon}\bigr)^{k-i+1} \bigl(\eb (1+
\eta_{i-1})^{\a+\epsilon}-\eb\bigl(\eta_{i-1}^{\a+\epsilon}
\bigr) \bigr)
\nonumber\\[-8pt]\\[-8pt]
&\leq& c\sum_{i=1}^k\bigl(\eb
A^{\a+\epsilon}\bigr)^{k-i+1} \leq c \frac{(\eb A^{\a+\epsilon})^{k}}{\eb A^{\a+\epsilon
}-1}.\nonumber
\end{eqnarray}
Indeed, we will prove that if $\epsilon< 1 $ then there is a constant $c$
such that for $i\ge1$,
\[
\eb(1+\eta_{i})^{\a+\epsilon}-\eb\eta_{i}^{\a+\epsilon}
\leq c.
\]
If $\a+\epsilon\leq1$ then
this follows from the concavity of the function $f(x)=x^{\a+\epsilon}$, $x>0$.
If $\a+\epsilon>1$ we
use the mean value theorem
to obtain
\[
\eb(1+\eta_{i})^{\a+\epsilon}-\eb\eta_{i}^{\a+\epsilon}
\leq(\a+\epsilon) \eb(1+\eta_{i})^{\a+\epsilon-1 }\leq(\a+\epsilon)\eb
\eta_\infty^{\a
+\epsilon-1 }<\infty.
\]
Now we choose $ \epsilon=k^{-0.5}$. Then by (\ref{eq6}), (\ref
{eq7}) and Lemma~\ref{exp},
\[
\eb\eta_k^{\a}\one_{ \{\eta_k> x\}}\leq c
\frac{(\eb A^{\a+\epsilon})^{k}}{\eb A^{\a+\epsilon}-1}x^{-\epsilon} \leq c \frac{ e^{\rho n_1/\sqrt{k}-\log
x/\sqrt{k}}}{\eb A^{\a+\epsilon}-1} \leq c
\frac{ e^{-\rho m/\sqrt{k}}}{\eb A^{\a+\epsilon}-1}.
\]
In the last step we used that
$k\leq n_1=n_0-m$, where $n_0=[\rho^{-1}\log x]$.
Moreover, since $m=[(\log x)^{0.5 +\sigma}]$,
$m/\sqrt{k}\geq2 c_1(\log x)^{\sigma}$ for some $c_1>0$.
On the other hand, setting $\gamma=\epsilon=k^{-0.5}$ in (\ref{eqtaylor}),
we obtain $\eb A^{\a+\epsilon}-1 \ge\rho k^{-0.5}/2$.
Combining the bounds above, we finally arrive
at
\[
{\sup_{r_n\le k\le n_1,b_n\le x}}k^{-1}\eb\eta_k^{\a}
\one_{\{\eta
_k> x\}}\leq c e^{-c_1 (\log x)^{\sigma}}
\]
for constants $c,c_1>0$.
This estimate yields the desired relation (\ref{eqnull}) and thus
completes the proof of the first part
of the lemma when $c_\infty^+ >0$.

If $c_\infty^+ = 0$ we proceed in the same way, observing that for any
$\d, z>0$ and sufficiently
large $x$,
\[
\frac{\pb\{Y > x/z\}}{ \P\{ |Y|>x \}} < \d z^\a
\]
and hence $I_1$ converges to 0 as $n$ goes to infinity.
This proves the lemma.
\end{appendix}

\section*{Acknowledgments}

We would like to thank the referee whose comments
helped us to improve the presentation of this paper.

%suskaldyti doi

% imsref loaded by lrinkeviciute, 2013-01-02 12:36:48

\printaddresses

\end{document}